\documentclass[aop,preprint]{imsart}
\RequirePackage{amsthm,amsmath,amsfonts,amssymb}
\RequirePackage[numbers,sort&compress]{natbib}
\RequirePackage[colorlinks,citecolor=blue,urlcolor=blue,linkcolor=blue]{hyperref}
\RequirePackage[mathscr]{eucal}
\usepackage{enumitem}
\RequirePackage{macro}
\numberwithin{equation}{section}

\allowdisplaybreaks

\begin{document}
\begin{frontmatter}
\title{Quantitative approximations for finite-state mean-field configuration models}
\runtitle{Quantitative mean-field approximations}
\begin{aug}
\author[A]{\fnms{Arnab}~\snm{Ganguly}\ead[label=e1]{aganguly@lsu.edu}}
\author[A]{\fnms{Yangrui}~\snm{Xiang}\ead[label=e2]{yxiang8@lsu.edu}}
\address[A]{Department of Mathematics, Louisiana State University\printead[presep={,\ }]{e1,e2}}
\end{aug}
\begin{abstract}
We establish quantitative approximations for the empirical measures of finite-state mean-field particle systems in a general framework. The scales of the empirical drift and covariance distinguish
the ODE and critical diffusion regimes and determine the
corresponding time rescaling.

In the ODE regime, we obtain quantitative
central limit approximations in the $1$-Wasserstein distance and a Bernstein-type concentration inequality. The concentration bound has no dimensional prefactor. Its explicit time
dependence also identifies growing intervals on which the
ODE approximation remains accurate.
In the critical regime, we prove a quantitative nonlinear diffusion approximation in the $1$-Wasserstein distance. 
We apply our quantitative approximation results to noisy voter systems, Moran models with mutations and selections, heat-bath Glauber dynamics for the $p$-spin Curie--Weiss model, and autocatalytic reactions with mutations.
\end{abstract}
\end{frontmatter}

\section{Introduction}\label{sec:intro}
\subsection{Motivation and related literature}
\label{subsec:intro_motivation_PRM}
We study interacting particle systems on configuration spaces, whose transition rates
depend on the states of other particles.
Classical references for interacting particle
systems and their scaling limits are 
Liggett~\cite{Liggett1985}, and Kipnis and
Landim~\cite{KipnisLandim1999}, respectively. Configurations record state of particle
at every site. 
Although each site has only finitely many possible states, the number of configurations can grow exponentially with the total number of particles. A central problem is therefore to obtain a
macroscopic description of the evolution of quantities such as particle densities.

This paper focuses on finite-state mean-field systems in which
each particle's transition rates depend on its current state and the empirical distribution of states. The empirical measure is then a finite-dimensional Markov process describing the proportions of particles in the different states. The primary goal of this article is to develop a general framework for quantitative central limit theorems, concentration inequalities without a dimensional prefactor, and quantitative nonlinear diffusion approximations. We use the drift and covariance scales to identify the relevant time rescaling and distinguish between regimes admitting ODE and SDE approximations. Importantly, our focus is on quantitative error estimates for
finite particle systems, rather than on large-$n$ asymptotics. Our framework covers several important models, including noisy voter systems, Moran models with mutations and selections, heat-bath Glauber dynamics for the $p$-spin Curie--Weiss model, and autocatalytic reactions with mutations.
These examples are discussed in Section~\ref{sec:examples_PRM}.

Classical hydrodynamic and mean-field limit theory provides
deterministic macroscopic descriptions through PDEs and ODEs,
respectively; see \cite{KipnisLandim1999,Kurtz1970,Kurtz1971}.
Quantitative estimates for fluctuations under equilibrium have
been studied in several related models.
For multivariate Markov jump processes near a stable
equilibrium, Barbour, Luczak and
Xia~\cite{BarbourLuczakXia2018I,BarbourLuczakXia2018II}
obtained approximation by discrete multivariate normal laws with total variation error of order $n^{-1/2}\log n$.
For Curie--Weiss models, Shao and
Zhang~\cite{ShaoZhang2019} proved Gaussian Berry--Esseen
bounds of order $n^{-1/2}$ in the high-temperature regime. At the critical temperature and zero external field,
Chatterjee and Shao~\cite[Theorem~2.1]{ChatterjeeShao2011}
obtained an $O(n^{-1/2})$ bound in Kolmogorov distance
between the law of the rescaled magnetization
$n^{1/4}m_n$ and the non-Gaussian law with density
proportional to $e^{-x^4/12}$. For non-equilibrium fluctuations, Jara and Menezes~\cite{JaraMenezes2018} used Yau's relative entropy method to show that
the non-equilibrium density fluctuations, for a weakly asymmetric exclusion process in dimensions $d\leq3$, converge to the solution of a linear stochastic heat equation.
In this article, we study
the fluctuations of the empirical vector around a time-dependent ODE trajectory. We obtain a $1$-Wasserstein bound
for the full fluctuation vector, uniformly over time marginals
on each fixed interval and over deterministic initial conditions. The regularity assumptions concern the empirical drift
and covariance. However, individual transition rates are not required to be smooth. The central limit approximation retains the finite-$n$
coefficients and initial empirical state, and requires no convergence of these two quantities.

Nonasymptotic concentration bounds have been established
for several related models as well. Darling and
Norris~\cite{DarlingNorris2008} controlled the time supremum
of deviations from ODE trajectories by exponential bounds. The coordinatewise formulation contains a prefactor proportional to the dimension.
Del Moral and Rio~\cite{DelMoralRio2011} established Hoeffding-, Bennett, and
Bernstein-type concentration inequalities for bounded
test observables in discrete-generation mean-field
particle models.
For equilibrium Curie--Weiss configurations, Chatterjee and
Dey~\cite{ChatterjeeDey2010} obtained concentration
of the form
$\mathbb P(n^{1/4}|m_n|\geq u)\leq 2e^{-cu^4}$ at the critical inverse temperature and zero external field,
where $m_n$ denotes the magnetization. For canonical equilibrium configurations of symmetric
exclusion, Gaussian concentration follows from
Bobkov's inequality on slices of the discrete
cube~\cite[Theorem~2.1]{Bobkov2004}.
More recently, Hermon and Salez~\cite{HermonSalez2023}
established modified log-Sobolev inequalities
for configuration measures satisfying the stochastic
covering property, which yields Gaussian concentration
for Lipschitz observables. In our work, our concentration estimate controls the time supremum deviation of the empirical vector. It requires no equilibrium, contractivity, or
negative-dependence assumption. The concentration estimate also applies on growing intervals
$[0,T_n]$, with $T_n\to\infty$. Its explicit formula shows the time dependence, which
gives sufficient conditions on the growth of $T_n$, under which the empirical measure remains uniformly close to its ODE trajectory with high probability.

The macroscopic evolution need not be deterministic.
Surprisingly, under suitable time rescaling, the empirical density itself can converge to a nondegenerate diffusion, so that an SDE describes its evolution. There are examples include the Wright--Fisher
limits for voter-model densities established by
Cox~\cite[Theorem~2]{Cox1989} and Chen, Choi and
Cox~\cite{ChenChoiCox2016}. For the noisy voter model in the weak noise regime, Aljovin, Jara and Xiang
\cite[Theorem~3.1]{AljovinJaraXiang2024} obtained a quantitative Wright--Fisher approximation on an accelerated time scale $nt$.
Therefore, both deterministic and stochastic macroscopic descriptions
arise in the study of configuration dynamics. These results concern voter and noisy voter models, while in this article we consider
a general class of finite-state mean-field systems and obtain quantitative approximations.
The transition rates, parameter scales, and observation times determine the approximation regimes.  Our framework gives both ODE and SDE approximations by allowing
the drift and covariance to have separate deterministic scales.
These scales identify the approximation regime and the time rescaling order. In the critical regime, we
obtain quantitative Wasserstein bounds for approximation of the empirical measure itself by a nonlinear diffusion.

\subsection{Main results: scaling regimes and quantitative bounds}
\label{subsec:intro_scaling_results_PRM}

We consider $n$ particles with states in $\mathbb S=\{1,\ldots,K\}$. The configurations are $x=(x_1,\ldots,x_n)$ on the state space $\mathbb S^n=\{1,\ldots,K\}^n$.
Write $r_k=n^{-1}\sum_{i=1}^n\mathbf1_{\{x_i=k\}}$ for the density of particles in state $k$. If $x_i=k$, we write $H_{k,l}^{(n)}(r)$ as the transition rate of the particle $i$ changing from $k$ to $l\ne k$. The transition rates depend on the density of all types of particles.

We denote $\mathsf e_i$ as the $i$th coordinate vector in $\mathbb R^K$.
This update changes the empirical distribution from $r$ to
$r+\bar{\canvec}_{k,l}/n$, where
$\bar{\canvec}_{k,l}=\canvec_l-\canvec_k$.
 Consequently,
\begin{align*}
 F^{(n)}(r)=\sum_{k\ne l}r_kH_{k,l}^{(n)}(r)\bar{\canvec}_{k,l}, \qquad
 G^{(n)}(r)=\sum_{k\ne l}r_kH_{k,l}^{(n)}(r)
              \bar{\canvec}_{k,l}\bar{\canvec}_{k,l}^{\top}.
\end{align*}
The vector $F^{(n)}$ describes the expected net redistribution
of density among states, while $G^{(n)}/n$ describes
the covariance of the random jumps in the empirical
distribution.

We use positive sequences $\alpha(n)$ and $\beta(n)$ to describe the growth or decay of $F^{(n)}$ and $G^{(n)}$ as the total population $n$ increases. We use $\alpha(n)$ and $\beta(n)$ to normalize the drift and
covariance coefficients, respectively, in 
Assumption~\ref{assum:scaling-1}. We study the empirical measure on the time scale $t/\alpha(n)$.
On the time scale $t/\alpha(n)$, the drift of the empirical measure is uniformly bounded. On any fixed interval $[0,T]$, the covariance of the random jumps accumulates to a matrix
with trace bounded by $C_T\,\beta(n)/(n\alpha(n))$, where $C_T$ is independent of $n$.
We take $a_n^2=n\alpha(n)/\beta(n)$, and the limit of $a_n$ determines which approximation regime the process falls into.
When $a_n\to\infty$, the noise vanishes and the empirical measure can be approximated by an ODE. When $a_n^{-2}$ converges to a strictly positive
finite constant, we obtain
a nonlinear diffusion approximation. We refer to the two cases above as
the ODE and critical diffusion regimes, respectively.

In the ODE regime, $\rho^{(n)}$ is the the trajectory starting from the same initial condition obtained from the normalized drift. Theorem~\ref{th:ODE-approx} gives the bound of the $p$th moment of the largest distance over a fixed time interval between the empirical measure and the ODE trajectory. 
The next central limit approximation describes the
fluctuations around $\rho^{(n)}$. We multiply by the proper scale $a_n$. The fluctuations
are approximated by a linear Gaussian diffusion.
Theorem~\ref{thm:kantorovich_PRM} bounds the $1$-Wasserstein distance between their time marginals by
$C_T\max\{a_n^{-1/2},(a_n/n)^{1/3}\}$,
uniformly for $0\leq t\leq T$. The constant $C_T$ is independent of $n$ and the initial empirical states. 

Theorem~\ref{thm:concentration_PRM} gives a Bernstein-type concentration inequality. The exponential bound gives the time supremum deviation of the empirical measure from its ODE trajectory, on a time interval $[0,T]$, for any $T>0$. 
The exponential bound depends explicitly on the covariance of the random jumps, the jump size, and the drift Lipschitz constant. Importantly, the exponential bound has no explicit dimensional factor. Consequently, the exponential bound is uniform
in $K$ whenever the normalized drift Lipschitz constant and covariance-trace bound do not depend on $K$.
The estimate also holds on growing intervals $[0,T_n]$,
with $T_n\to\infty$. The exponential bound gives sufficient conditions on the growth of $T_n$. Under such conditions, the exponential bound stays small enough for sufficiently large $n$. Therefore, with high probability, the empirical measure stays close to its ODE trajectory throughout $[0,T_n]$.

In the critical diffusion regime, we
approximate the empirical measure itself by a nonlinear diffusion. Under Assumption~\ref{assum:critical_regime_PRM},
Theorem~\ref{thm:critical_dK_PRM} bounds the $1$-Wasserstein distance between their time marginals by
$C_{\mathrm{SDE},1}(T)
 \bigl[n^{-1/3}+|a_n^{-1}-a_\infty^{-1}|\bigr]$,
uniformly over $[0,T]$. The constant $C_{\mathrm{SDE},1}(T)$ depends on $a_*,T,K$, and the relevant drift and covariance derivative bounds.
In Theorem \ref{thm:critical_pairwise}, the additional regularity and balance conditions 
can improve the quantitative bound from $n^{-1/3}$ to $n^{-1/2}$. For Wright--Fisher-type diffusions, the smoothness assumption on
the diffusion coefficient fails at the boundary, However, we can still give a quantitative approximation bound for such diffusions.
Theorem~\ref{thm:critical_WF_dK_PRM} gives the upper bound
$C_{\mathrm{WF},1}(T)[n^{-1/3}+\varepsilon_n^{1/2}
+|a_n^{-2}-a_\infty^{-2}|^{1/2}]$ in $1$-Wasserstein distance,
uniformly over $[0,T]$, where $\varepsilon_n$ is the covariance
discrepancy in \eqref{eq:WF_covariance_error_PRM} and
$n^{-1/3}$ can be improved to $n^{-1/2}$ under
\eqref{eq:pairwise_balance}.

Theorems~\ref{th:ODE-approx}, \ref{thm:kantorovich_PRM},
\ref{thm:concentration_PRM}, \ref{thm:critical_dK_PRM}, \ref{thm:critical_pairwise}, and \ref{thm:critical_WF_dK_PRM},
allow the coefficients and initial empirical states to depend
on $n$. Their convergence is not required. The proofs use the coefficients and initial
state of each finite-$n$ system.
Therefore, all quantitative bounds in this article hold for all $n$ sufficiently large, without passing to the limit when $n \to \infty$.

\subsection{Applications of the main theorems}
\label{subsec:intro_methods_applications_PRM}
The examples in Section~\ref{sec:examples_PRM} show how parameter scales determine the approximation regime in four classes of
interacting particle systems. 

Let $(\gamma_n)_{n\geq1}$ be a positive sequence.
In the noisy voter model, the noise given by the spontaneous flip, has rates proportional
to $\gamma_n/n$. In the Moran model with selections, mutation rates and the increases in reproduction rates promoted by advantageous selections are proportional to $\gamma_n/n$.
If $\gamma_n\to\infty$ and $\sup_n(\gamma_n/n)<\infty$, the ODE-regime theorems apply on the time scale $t\mapsto nt/\gamma_n$.
If instead $\gamma_n\to\gamma_\infty\in(0,\infty)$, the time
scale $t\mapsto nt$ gives the critical diffusion regime, with quantitative approximations by a Wright--Fisher diffusion.

For heat-bath Glauber dynamics of the $p$-spin Curie--Weiss model, the system is in the ODE
regime. The ODE regimes theorems apply on
each fixed interval of the original time scale. This covers values of the inverse temperature and external magnetic field $(\beta,h)$ at which the equilibrium magnetization undergoes a phase transition.

For the autocatalytic reactions with mutations, per-individual conversion rates are of constant order
and mutation rates are of order $n^{-1}$.
In the balanced conversion case, the system is in the critical diffusion regime on the time
scale $t\mapsto nt$. For the unbalanced conversion case, the system falls into the ODE regime on the
original time scale.

\subsection{Organization of the paper}
\label{subsec:intro_organization_PRM}
Section~\ref{sec:framework_main_results_PRM} states the main results.
Section~\ref{sec:ODE_regime_PRM} gives the proofs in the ODE regime,
and Section~\ref{sec:critical_nonlinear_diffusion_PRM} gives the proofs
in the critical SDE regime. Section~\ref{sec:examples_PRM} applies our main results to four classes of models: noisy voter systems, Moran models with mutation and selection, heat-bath Glauber dynamics for the $p$-spin
Curie--Weiss model, and autocatalytic reactions with mutations. Appendix~\ref{app:regularization} contains
the regularization lemma for Wasserstein distance and proofs for regularity lemmas in Section~\ref{sec:ODE_regime_PRM} and Section~\ref{sec:critical_nonlinear_diffusion_PRM}.

\textbf{Notations:} 
We write $\|\cdot\|$ for the Euclidean norm. For an $m$-linear
map $A:(\R^d)^m\to\R^q$, its operator norm is defined by
$$
\|A\|_{\mathrm{op}}
\dfeq
\sup_{\|v_1\|\leq1,\ldots,\|v_m\|\leq1}
\|A[v_1,\ldots,v_m]\|.
$$
For scalar-valued maps, the norm on the right is the absolute value.
For a matrix $A\in\R^{q\times d}$, this gives the usual matrix
operator norm,
$$
\|A\|_{\mathrm{op}}
=\sup_{\|v\|\leq1}\|Av\|
=\sup_{\|u\|\leq1,\|v\|\leq1}|u^\top Av|,
$$
where $u\in\R^q$ and $v\in\R^d$. The trace (or nuclear) norm
of $A$ is defined by
$$
\|A\|_{\mathrm{tr}}
\dfeq \trace\bigl((A^\top A)^{1/2}\bigr).
$$

For an integer $m\geq0$,  $C^m(\R^K;\R^q)$ will denote the space
of $m$ times continuously differentiable maps from $\R^K$ to $\R^q$.
For $f\in C^m(\R^K;\R^q)$ and $1\leq j\leq m$, we regard
$D^j f(x)$ as a symmetric $j$-linear map from $(\R^K)^j$
to $\R^q$ and, unless a different domain of the supremum is specified, write
$$
\|D^j f\|_{\mathrm{op},\infty}
\dfeq
\sup_{x\in\R^K}\|D^j f(x)\|_{\mathrm{op}},
\qquad
\|f\|_\infty
\dfeq
\sup_{x\in\R^K}\|f(x)\|.
$$
We identify $Df(x)$ with its Jacobian matrix, so that
$\|Df(x)\|_{\mathrm{op}}$ is the usual matrix operator norm. When $q=1$, $Df(x)$ is identified with $\nabla f(x)$ and the operator norm is simply the Euclidean norm, and we  write $\|Df(x)\|_{\mathrm{op}} \equiv \|Df(x)\| =\|\nabla f(x)\|$. 
We consider the space $C_b^m(\R^K;\R^q)$ consists of those
$f\in C^m(\R^K;\R^q)$ for which
$$
\|f\|_{C_b^m}
\dfeq
\|f\|_\infty
+\sum_{j=1}^m\|D^j f\|_{\mathrm{op},\infty}
<\infty.
$$

For probability measures on $\mathbb R^K$ with finite first moments,
the \emph{Wasserstein-1 distance} (or Kantorovich distance) is
\begin{equation}\label{eq:dK_definition}
\wass_1(\nu,\widetilde\nu)
 =\sup_{\operatorname{Lip}(f)\leq1}
   \left|\int f\,d\nu-\int f\,d\widetilde\nu\right|,
\end{equation}
where Lipschitz constants are taken in the Euclidean norm.

\section{Framework and main results}
\label{sec:framework_main_results_PRM}
Throughout this article, $K\geq2$ is fixed. We set that
\[
 \mathbb S=\{1,\ldots,K\},
 \qquad
 \mathcal P(\mathbb S)
 =\left\{r\in[0,1]^K:\sum_{k=1}^K r_k=1\right\}.
\]

\subsection{Configuration process and empirical characteristics}
We first define the empirical measure map
$\nden:\mathbb S^n\to\mathcal P(\mathbb S)$ as
\begin{equation}\label{eq:den}
\nden(x)=\frac1n\sum_{i=1}^n\delta_{x_i},
\qquad x\in\mathbb S^n.
\end{equation}
Let $\nH_{k,l}:\mathcal P(\mathbb S)\to[0,\infty)$, $k\ne l$,
be finite measurable transition-rate functions. For $1\le i\le n$
and $k\ne l$, define the jump maps
$\phi^{(i,k,l)}:\mathbb S^n\to\mathbb S^n$ by
$$
\phi^{(i,k,l)}(x)
=x+(l-k)\mathsf e_i\mathbf1_{\{x_i=k\}},
$$
where $\mathsf e_i$ is the $i$th standard basis vector of
$\mathbb R^n$.

We consider a continuous-time Markov chain (CTMC)
$X^{(n)}=(X_1^{(n)},\ldots,X_n^{(n)})$ on $\mathbb S^n$
with transitions
$$
x\longrightarrow\phi^{(i,k,l)}(x)
\quad\text{at rate}\quad
\nH_{k,l}(\nden(x)),
\qquad x_i=k,\quad l\ne k.
$$
We now give a pathwise representation of $X^{(n)}$.  Let $\xi_{i,k,l}$ be independent Poisson random measures (PRMs) on $[0,\infty)^2$ with intensity measure $\leb \otimes \leb$, independent of $\nX(0)$, and let
$\widetilde\xi_{i,k,l}=\xi_{i,k,l}-\leb \otimes \leb$ denote the corresponding compensated PRMs. Then $X^{(n)}$ satisfies the following SDE system:
\begin{equation}\label{eq:SDE-rep}
 \nX_i(t)=\nX_i(0)+\sum_{k\ne l}(l-k)
 \int_0^t\int_0^\infty
 \mathbf1_{\{u\leq\nH_{k,l}(\nden(\nX(s-)))\}}
 \mathbf1_{\{\nX_i(s-)=k\}}\xi_{i,k,l}(du\,ds).
\end{equation}
The generator of the Markov process $X^{(n)}$ is given by
\begin{equation}\label{eq:X-gen}
 \ngen f(x)=\sum_{i=1}^n\sum_{k\ne l}
 \nH_{k,l}(\nden(x))\mathbf1_{\{x_i=k\}}
 \bigl[f(\phi^{(i,k,l)}(x))-f(x)\bigr].
\end{equation}

For $f:\mathbb S\to\mathbb R$, 
\begin{equation}\label{eq:Ito-i}
 f(\nX_i(t))=f(\nX_i(0))+\int_0^tL_i^{(n)}f(\nX(s))\,ds+\nM_i(f,t),
\end{equation}
where the operator $L_i^{(n)}$ is given by
\begin{equation}\label{eq:gen-i}
 L_i^{(n)}f(x)=\sum_{k\ne l}\nH_{k,l}(\nden(x))
 [f(l)-f(k)]\mathbf1_{\{x_i=k\}},
\end{equation}
and  the martingale $\nM_i(f,\cdot)$ is defined by
\begin{equation}
\label{eq:mart-i}
\begin{aligned}
 \nM_i(f,t)
 ={}&\sum_{\substack{k,l\\k\neq l}}^K
       \bigl(f(l)-f(k)\bigr)
 \int_{[0,\infty)\times[0,t]}
 \indic_{[0,\nH_{k,l}(\nden(\nX(s-)))]}(u) \\
 & \hspace{2cm} \times
 \indic_{\{\nX_i(s-)=k\}}\,
 \widetilde\xi_{i,k,l}(du\times ds).
\end{aligned}
\end{equation}
Consequently, summing \eqref{eq:Ito-i} gives the following equation for  $\mu^{(n)}_\cdot \dfeq \nden(\nX(\cdot))$, the empirical measure of the particle-system $\nX$:
\begin{align} \label{eq:char-emp-meas}
\mu_t^{(n)}(f)=&\mu_0^{(n)}(f)
 +\int_0^t F^{(n)}(f,\mu_s^{(n)})\,ds+\nmart(f,t),
\end{align}
where for $f: \mathbb{S} \rt \R$, the function $F^{(n)}(f,\cdot):  \mathcal P(\mathbb{S}) \rt \R$ is defined as 
\begin{align} \label{eq:drft-def}
F^{(n)}(f, r) \dfeq  \sum_{ \substack{k, l \\ k \neq l}}^K (f(l) - f(k)) r_k \nH_{k,l}(r), \quad r=(r_1,r_2,\hdots, r_K) \in \mathcal P(\mathbb{S}),
\end{align}
and the martingale $ \nmart(f,\cdot)$ is given by
\begin{equation}
\label{eq:empirical_martingale_PRM}
\begin{aligned}
 \nmart(f,t)
 \equiv \frac1n\sum_{i=1}^n\nM_i(f,t)
 =&\frac1n\sum_{i=1}^n
       \sum_{\substack{k,l\\k\neq l}}^K
       \int_{[0,\infty)\times[0,t]}
       [f(l)-f(k)]\\
 &\qquad\times
       \mathbf1_{\{\nX_i(s-)=k\}}
       \mathbf1_{\{u\leq\nH_{k,l}(\mu_{s-}^{(n)})\}}\,
       \widetilde\xi_{i,k,l}(du,ds).
\end{aligned}
\end{equation}
Notice that the predictable  covariation between $\nmart(f,\cdot)$ and $\nmart(g,\cdot)$ is given by
\begin{align} \label{eq:mart-emp-quad}
\begin{aligned}
    \langle\nmart(f,\cdot),\nmart(g,\cdot)\rangle_t
 =&\ \f{1}{n} \sum_{ \substack{k, l \\ k \neq l}}^K \big(f(l) - f(k)\big)\big(g(l) - g(k)\big) \int_0^t \nH_{k,l}(\mu^{(n)}_s) \mu^{(n)}_s(k)\ ds\\
\dfeq&  \f{1}{n}\int_0^tG^{(n)}(f,g, \mu^{(n)}_s) ds
\end{aligned}
\end{align}
Equivalently, identifying $\mu^{(n)}$ as a vector-valued process with state space $\mathcal{P}(\mathbb{S}) \subset \R^K$,  \eqref{eq:char-emp-meas}--\eqref{eq:mart-emp-quad}, in vector notations become
\begin{equation}\label{eq:char-emp-meas-scaled-vec}\mu_t^{(n)}=\mu_0^{(n)}+\int_0^tF^{(n)}(\mu_s^{(n)})\,ds
 +\mart^{(n)}(t),
\end{equation}
where $F^{(n)}(r) \dfeq \lf(F^{(n)}(k,r)\equiv F^{(n)}(\indic_{\{k\}},r): k=1,2\hdots K\ri) \in \R^K$ and

$G^{(n)}(r) \dfeq \lf(G^{(n)}(k,l,r)=G^{(n)}(\mathbf1_{\{k\}},\mathbf1_{\{l\}},r): 1\leq k,l \leq K\ri)$ are given by
\begin{equation}\label{eq:cov-def-comp}
F^{(n)}(r)=\sum_{ \substack{k, l \\ k \neq l}}r_k\nH_{k,l}(r)\bar{\canvec}_{k,l},\qquad
 G^{(n)}(r)=\sum_{ \substack{k, l \\ k \neq l}}r_k\nH_{k,l}(r)\bar{\canvec}_{k,l}\bar{\canvec}_{k,l}^\top, \qquad \bar{\canvec}_{k,l} \dfeq \canvec_l - \canvec_k
\end{equation}
and the martingale $\nmart$ is given by
\begin{equation}\label{eq:empirical_martingale_PRM-2}
 \nmart(t)\equiv \frac{1}{n}\sum_{i=1}^n\sum_{ \substack{k, l \\ k \neq l}}^K
 \int_{[0,\infty)\times[0,t]} \bar{\canvec}_{k,l}
 \mathbf1_{\{\nX_i(s-)=k\}}
 \mathbf1_{\{u\leq\nH_{k,l}(\mu_{s-}^{(n)})\}}
 \,\widetilde\xi_{i,k,l}(du,ds),
\end{equation}
with the matrix covariation given by 
$$\langle\mart^{(n)}\rangle_t=\frac1n\int_0^tG^{(n)}(\mu_s^{(n)})\,ds.$$
We also write $F_k^{(n)}(r)=F^{(n)}(k,r)$.
Notice that the trace of the nonnegative definite (n.n.d.) matrix $G^{(n)}$ is given by
\begin{equation}\label{eq:trace_total_rate}
\trace G^{(n)}(r)=\sum_{\substack{k,l=1\\k\neq l}}^K
   r_k\nH_{k,l}(r)\lf\|\bar{\canvec}_{k,l}\ri\|^2 = 2\sum_{ \substack{k, l=1 \\ k \neq l}}^Kr_k\nH_{k,l}(r).
\end{equation}

\subsection{Scaling regimes}
For simplicity, we assume throughout that the initial empirical measure $\mu_0^{(n)}$ is deterministic (see Remark \ref{rem:random_initial_PRM}). 
The scaling regimes considered in this paper depend on the growth rates
$\alpha(n)$ and $\beta(n)$ of the drift term $F^{(n)}$ and the covariance
term $G^{(n)}$, respectively. We impose the following assumption
throughout the paper.

\begin{assumption}[General scaling]
\label{assum:scaling-1}
There exist positive sequences $\alpha(n)$ and $\beta(n)$ such that
the normalized coefficients
\begin{equation}\label{eq:scaled-drft-quad-def-1}
 \bar F^{(n)}\dfeq (\alpha(n))^{-1} F^{(n)},\qquad
 \bar G^{(n)}\dfeq (\beta(n))^{-1}G^{(n)}
\end{equation}
satisfy
\begin{equation}\label{eq:scaled-drft-quad-bdd}
\|\bar F\|_{\infty} \dfeq \sup_{n\geq1}\sup_{r\in\mathcal P(\mathbb S)}
 \|\bar F^{(n)}(r)\|< \infty, \quad \|\bar G\|_{\mathrm{tr},\infty} \dfeq \sup_{n\geq1}\sup_{r\in\mathcal P(\mathbb S)}\trace(\bar G^{(n)}(r))<\infty.
\end{equation}
\end{assumption}
Let $\bnX(\cdot) \dfeq \nX(\cdot/\alpha(n))$.  The main object of our study is the time-scaled empirical measure
\begin{align}
    \label{eq:emp-meas-scld}
    \bnmu_\cdot \dfeq \mu_{\cdot/\alpha(n)}^{(n)} = \nden(\bnX(\cdot)).
\end{align}
Notice that a simple time change in \eqref{eq:SDE-rep} gives
the following SDE for $\bnX=(\bnX_1,\ldots,\bnX_n)$:
\begin{align}
\label{eq:SDE-rep-scaled}
\begin{aligned}
\bnX_i(t) =& \nX_i(0)+ \sum_{ \substack{k, l \\ k \neq l}}^K \int_{[0,\infty)\times[0,t]} (l-k) \indic_{[0,\nH_{k,l}(\nden(\bnX(s-)))]}(u)\\
& \hspace{1.5cm}\times \indic_{\{\bnX_i(s-) = k\}} \xi_{i,k,l}\lf(du\times \alpha(n)^{-1}ds\ri),
\end{aligned}
\end{align}
with  the generator $\bngen$ of $\bnX$ satisfying $\bngen \equiv \f{1}{\alpha(n)} \ngen$. For $f:\R^K\to\R$, setting $r=\nden(x)$,  we have  
\begin{align}\label{eq:scaled_empirical_generator_PRM}
\bngen(f\circ\nden)(x)
&=\frac1{\alpha(n)}
  \sum_{i=1}^n\sum_{\substack{k,l=1\\k\neq l}}^K
  \indic_{\{x_i=k\}}\nH_{k,l}(r)
  \left[
    f\left(r+\frac{\bar{\canvec}_{k,l}}n\right)-f(r)
  \right]\notag\\
&=\frac n{\alpha(n)}
  \sum_{\substack{k,l=1\\k\neq l}}^K
  r_k\nH_{k,l}(r)
  \left[
    f\left(r+\frac{\bar{\canvec}_{k,l}}n\right)-f(r)
  \right].
\end{align}
\eqref{eq:SDE-rep-scaled} shows that $\bnmu$ satisfies
\begin{equation}\label{eq:scaled_characteristics_main}
\bnmu_t
=\mu_0^{(n)}
+\int_0^t\bar F^{(n)}(\bnmu_s)\,ds
+\bnmart(t),
\end{equation}
where the martingale $\bnmart(\cdot)\dfeq \nmart(\cdot/\alpha(n))$ admits the expression
\begin{align}\label{eq:scaled_empirical_martingale_PRM}
\begin{aligned}
\bnmart(t)
&=\frac1n\sum_{i=1}^n
\sum_{\substack{k,l=1\\k\neq l}}^K
\int_{[0,\infty)\times[0,t]}
\bar{\canvec}_{k,l}\,
\indic_{\{\bnX_i(s-)=k\}}
\indic_{[0,\nH_{k,l}(\bnmu_{s-})]}(u)\,
\widetilde\xi_{i,k,l}
\lf(du\times\alpha(n)^{-1}ds\ri),
\end{aligned}
\end{align}
and has the predictable covariation matrix given by
\begin{equation}\label{eq:scaled_empirical_covariation}
\langle\bnmart\rangle_t
=\frac1{n\alpha(n)}
\int_0^tG^{(n)}(\bnmu_s)\,ds
=\frac{\beta(n)}{n\alpha(n)}
\int_0^t\bar G^{(n)}(\bnmu_s)\,ds.
\end{equation}
Here
$
\widetilde\xi_{i,k,l}
\lf(du\times\alpha(n)^{-1}ds\ri)
=
\xi_{i,k,l}
\lf(du\times\alpha(n)^{-1}ds\ri)
-\alpha(n)^{-1}du\,ds.
$
The appropriate approximation for $\bnmu$ depends on the scaling sequence
\begin{align}
    \label{eq:ef_scale}
    a_n \dfeq \left(\frac{n\alpha(n)}{\beta(n)}\right)^{1/2}.
\end{align}
Specifically, we obtain nonasymptotic bounds for approximating
$\bnmu$ by an ODE with vector field $\bnF$ when $a_n\to\infty$,
and by an SDE with drift term $\bnF$ and diffusion coefficient
$a_n^{-1}(\bnG)^{1/2}$ when $a_n\to a_\infty \in(0,\infty)$. 

\subsection{The ODE regime: \texorpdfstring{$a_n \to \infty$}{an to infinity}}\label{ode:main results}
\begin{assumption}[ODE regime]
\label{assum:ode-regime}
The following conditions hold.
\begin{enumerate}[label=(\roman*), ref=(\roman*)]
 \item\label{item:an-infty}
    $a_n \to \infty$ as $n \to \infty$;
  \item \label{item:barF-Lipschitz}  There exists an open set
$\mathbb U\subset\R^K$, independent of $n$, containing
$\mathcal P(\mathbb S)$, such that, for every $n\geq1$,
there exists an extension $\bnF\in C^1(\mathbb U;\R^K)$ satisfying
$$
\left.\bnF\right|_{\mathcal P(\mathbb S)}
=\alpha(n)^{-1}F^{(n)}.
$$
Furthermore, there exist finite constants $L_{DF}, L_G \geq 0$,
independent of $n$, such that
\begin{align}\label{eq:barF_G_reg_PRM}
\begin{aligned}
&\|D\bar F\|_{\mathrm{op},\infty}
\dfeq
\sup_{n\geq1}\sup_{r\in\mathcal P(\mathbb S)}
\|D\bnF(r)\|_{\mathrm{op}}<\infty,\\
&\|D\bnF(r)-D\bnF(r')\|_{\mathrm{op}}
\ \leq\ L_{DF}\|r-r'\|,\\
&\|\bnG(r)-\bnG(r')\|_{\mathrm{tr}}
\ \leq\ L_G\|r-r'\|,
\end{aligned}
\end{align}
for all $n\geq1$ and $r,r'\in\mathcal P(\mathbb S)$.
\end{enumerate}
\end{assumption}

Let $\rho^{(n)}$ be the solution to the ODE
\begin{equation}\label{eq:ODE-approx}
 \dot\rho_t^{(n)}=\bar F^{(n)}(\rho_t^{(n)}),\qquad
 \rho_0^{(n)}=\mu_0^{(n)}.
\end{equation}
This solution is unique and remains in $\mathcal P(\mathbb S)$: the
drift preserves total mass and its $k$th component is nonnegative on
the face $r_k=0$. A direct justification is given in
Section~\ref{sec:ode_moment_proof}.

\begin{theorem}[Quantitative ODE approximation]\label{th:ODE-approx}
Suppose that Assumption~\ref{assum:scaling-1} and Assumption~\ref{assum:ode-regime} hold. Then the following hold:
\begin{enumerate}[label=(\alph*), ref=(\alph*)]
\item \label{item:limit-is-prob} $\rho_t^{(n)} \in \mathcal P(\mathbb S)$ for all $t \in [0,T]$.
\item \label{item:ODE-approx} For any $p \geq 2$, there exists a constant $C_{\mathrm{ODE}}(p,T) \equiv C_{\mathrm{ODE}}(p,T, \|D\bar F\|_{\mathrm{op},\infty}, \|\bar G\|_{\mathrm{tr},\infty} )  > 0$ such that
\begin{equation}\label{eq:ODE-error-est}
 \mathbb E\left[\sup_{t\leq T}\|\bnmu_t-\rho_t^{(n)}\|^p\right]
 \leq C_{\mathrm{ODE}}(p,T)(a_n^{-p}+n^{-p}) \leq\ 2C_{\mathrm{ODE}}(p,T)(a_n \wedge n)^{-p}
\end{equation}
\end{enumerate}
\end{theorem}
Theorem~\ref{th:ODE-approx} is proved in Section~\ref{sec:ode_moment_proof}.

We now prove a concentration inequality for the ODE approximation.
The bound reflects two scales: $a_n^{-2}$ for the predictable quadratic
covariation of the martingale $\bnmart$ on fixed time intervals,
and $n^{-1}$ for its jump sizes.
\begin{theorem}[Exponential concentration]
\label{thm:concentration_PRM}
Under Assumption~\ref{assum:scaling-1} and Assumption~\ref{assum:ode-regime}:\ref{item:an-infty} \& \ref{item:barF-Lipschitz}, for every $n\geq1$ and $T,\delta>0$,
\begin{equation}\label{eq:concentration_main_PRM}
\begin{aligned}
 \mathbb P\left(
 \sup_{t\leq T}\|\bnmu_t-\rho_t^{(n)}\|\geq\delta
 \right)
 \leq&\ 2\exp\left\{
 -\frac{\delta_T^2a_n^2}
 {2\bigl(\|\bar G\|_{\mathrm{tr},\infty}T+\sqrt2\,\delta_Ta_n^2/(3n)\bigr)}
 \right\}\\
 \leq&\ 2\exp\{-C_{\mathrm{conc}}(T,\delta)(a_n^2\wedge n)\},
 \end{aligned}
\end{equation}
where 
\[
\delta_T=\delta e^{-\|D\bar F\|_{\mathrm{op},\infty}T}, \qquad  C_{\mathrm{conc}}(T,\delta)
  \dfeq \frac{\delta^2_T}
 {2\bigl(\|\bar G\|_{\mathrm{tr},\infty}T+\sqrt2\,\delta_T/3\bigr)}.
\]
\end{theorem}

Theorem~\ref{thm:concentration_PRM} is proved in
Section~\ref{sec:concentration_proof}.

Our final result in this regime is the central limit approximation of $\bnmu$. 
Define the fluctuation process by
\begin{equation}\label{eq:xi_PRM_def}
 \fluct_t^{(n)}=a_n(\bnmu_t-\rho_t^{(n)}).
\end{equation}

 Let $W$ be a standard $K$-dimensional Brownian motion, and let $Y^{(n)}_{s,x}(t)$  the solution to the time-inhomogeneous linear SDE
\begin{equation}\label{eq:Y_PRM_SDE}
 dY^{(n)}_{s,x}(t)=D\bar F^{(n)}(\rho_t^{(n)})Y^{(n)}_{s,x}(t)\,dt
       +(\bar G^{(n)})^{1/2}(\rho_t^{(n)})\,dW(t),\quad s\leq t, \quad Y^{(n)}_{s,x}(s)=x,
\end{equation}
with the (time-dependent) generator 
\begin{equation}\label{eq:Lscript_PRM}
 \mathcal{A}_t^{(n)}\psi(x)=D\bar F^{(n)}(\rho_t^{(n)})x\cdot D\psi(x)
  +\frac{1}{2} \trace\lf(\bnG(\rho_t^{(n)})D^2\psi(x)\ri).
\end{equation}
Here for $r \in \mathcal P(\mathbb S),$ $(\bar G^{(n)})^{1/2}(r)$ denotes the unique n.n.d. square root of the n.n.d. matrix $\bar G^{(n)}(r).$

\begin{theorem}[Central limit approximation in Wasserstein metric]
\label{thm:kantorovich_PRM}
Suppose that Assumption~\ref{assum:scaling-1} and Assumption~\ref{assum:ode-regime} hold.

Then there exists a finite constant $C_{\mathrm{CLT}}(T)\equiv C_{\mathrm{CLT}}\bigl(T, K,
 \|D\bar F\|_{\mathrm{op},\infty}, \|\bar G\|_{\mathrm{tr},\infty},
 L_{DF},L_G \bigr)$
(independent of $n$) such that
\begin{equation}\label{eq:dK_bound_PRM}
 \sup_{t\leq T} \wass_1\bigl(\operatorname{Law}(\fluct_t^{(n)}),
                     \operatorname{Law}(Y_{0,0}^{(n)}(t))\bigr)
 \leq C_{\mathrm{CLT}}(T)\max\left\{ a_n^{-1},a_n^{-1/2},(a_n/n)^{1/3} \right\}.
\end{equation}
\end{theorem}
Theorem~\ref{thm:kantorovich_PRM} is proved in Section~\ref{sec:gaussian_proof}.

\subsection{The critical diffusion regime}
Here the empirical measure itself is compared with a nonlinear diffusion.

\begin{assumption}[Critical regime]
\label{assum:critical_regime_PRM}
The following conditions hold.
\begin{enumerate}[label=(\roman*), ref=(\roman*)]
\item\label{item:an-a_infty} $a_* \dfeq \inf_{n\geq 1}a_n \in (0,\infty)$ and $\displaystyle a_n=\lf(\frac{n\alpha(n)}{\beta(n)}\ri)^{1/2}
 \longrightarrow a_\infty\in(0,\infty)$ as $n \to \infty.$
\item \label{item:FF-exten}
For every $n\geq1$, there is a factor
$\crdiff^{(n)}:\mathcal P(\mathbb S)\to\R^{K\times K}$ satisfying
$$
 \crdiff^{(n)}(r)\bigl(\crdiff^{(n)}(r)\bigr)^{\mathsf T}
 =\bnG(r),\qquad r\in\mathcal P(\mathbb S).
$$
The functions $\bnF$ and $\crdiff^{(n)}$ admit extensions,
denoted by the same symbols, such that
$$
 \bnF\in C^3(\R^K;\R^K),
 \qquad
 \crdiff^{(n)}\in C^3_b(\R^K;\R^{K\times K}),
$$
with $\|\crdiff\|^2_{\mathrm{op},\infty}\dfeq  \sup_{n\geq1}\sup_{u \in \R^K}\|\crdiff^{(n)}(u)\|^2_{\mathrm{op}}< \infty$
and for $j=1,2,3$
\begin{align}\label{eq:critical_uniform_C3_PRM}
\begin{aligned}
\|D^{j}\bar{F}\|_{\mathrm{op}, \infty} \dfeq&\  \sup_{n\geq1}\sup_{u \in \R^K}\|D^j\bnF(u)\|_{\mathrm{op}} < \infty, \\ \|D^{j}\Sigma\|_{\mathrm{op}, \infty} \dfeq&\  \sup_{n\geq1}\sup_{u \in \R^K}\|D^j\Sigma^{(n)}(u)\|_{\mathrm{op}} < \infty.
\end{aligned}
\end{align}
\end{enumerate}
\end{assumption}

We point out that, in the general critical-diffusion statements and proofs,
$\|D^j\bar F\|_{\mathrm{op},\infty}$ denotes the supremum
over $\R^K$ in \eqref{eq:critical_uniform_C3_PRM}.
In the ODE regime, $\|D\bar F\|_{\mathrm{op},\infty}$
retains the definition over $\mathcal P(\mathbb S)$ in
\eqref{eq:barF_G_reg_PRM}.

On $\R^K$, define
$$
 \bnG(u)\dfeq
 \crdiff^{(n)}(u)\bigl(\crdiff^{(n)}(u)\bigr)^{\mathsf T}.
$$
Note that $\Sigma^{(n)}$ need not be the symmetric n.n.d.
square root of $\bnG$, as a different factor may admit a
smooth extension even when $(\bnG(\cdot))^{1/2}$ does not.

\begin{remark}
Note that we can assume $\bnF,\crdiff^{(n)}\in C_b^3(\R^K)$
without loss of generality, with norms bounded uniformly in $n$.
This can be done by multiplying their $C^3$ extensions by a
fixed smooth cutoff which is one near $\mathcal P(\mathbb S)$
and has compact support in the domain of extension, and then
extending by zero. The same argument applies to the $C_b^4$
condition in  Theorem \ref{thm:critical_pairwise}.
\end{remark}

Let $W$ be a standard $K$-dimensional Brownian motion, and for $c>0$ let $\crlim_{z}^{(n,c)}$ be the solution to the SDE
\begin{equation}\label{eq:critical_diffusion_SDE_PRM}
 d\crlim_z^{(n,c)}(t)
 =\bnF(\crlim_z^{(n,c)}(t))\,dt
 +c^{-1}\crdiff^{(n)}(\crlim_z^{(n,c)}(t))\,dW(t),
 \qquad \crlim^{(n,c)}_z(0)=z,
\end{equation}
with the generator given by
\begin{align}\label{eq:critical_diffusion_generator_PRM}
 \mathcal B_c^{(n)}\psi(x)
 =\bnF(x)\cdot D\psi(x)
 +\frac{1}{2c^2}\trace\bigl(\bnG(x)D^2\psi(x)\bigr), \quad  \psi\in C^2(\R^K;\R)
\end{align}
Notice \eqref{eq:critical_diffusion_SDE_PRM} admits a unique strong solution in $\R^K$.

\begin{theorem}[Nonlinear diffusion approximation]
\label{thm:critical_dK_PRM}
Suppose that Assumption \ref{assum:scaling-1} and Assumption~\ref{assum:critical_regime_PRM} hold. Then, for every $T>0$, there exist finite constants
$C_{\mathrm{SDE},0}(T)$ and $C_{\mathrm{SDE},1}(T)$,
depending only on $a_*,T,K$,
$\|\crdiff\|_{\mathrm{op},\infty}$, and
$\|D^j\bar F\|_{\mathrm{op},\infty}$,
$\|D^j\crdiff\|_{\mathrm{op},\infty}$, $j=1,2,3$,
such that, for every $n\geq1$,
\begin{align}\label{eq:critical_dK_bound_PRM}
\begin{aligned}
 \sup_{t\leq T}\wass_1\lf(
 \operatorname{Law}(\bnmu_t),
 \operatorname{Law}\lf(\crlim^{(n,a_n)}_{\mu^{(n)}_0}(t)\ri)
 \ri) \leq & C_{\mathrm{SDE},0}(T)n^{-1/3}\\
 \sup_{t\leq T}
 \wass_1\lf(
 \operatorname{Law}(\bnmu_t),
 \operatorname{Law}\lf(\crlim^{(n,a_\infty)}_{\mu^{(n)}_0}(t)\ri)
 \ri)
 \leq&\ C_{\mathrm{SDE},1}(T)
 \bigl[n^{-1/3}+|a_n^{-1}-a_\infty^{-1}|\bigr].
\end{aligned}
\end{align}
\end{theorem}

Theorem~\ref{thm:critical_dK_PRM} is proved in
Section~\ref{sec:critical_comparison_proof}. The bounds in Theorem \ref{thm:critical_dK_PRM} can be
improved when the flows in opposite directions between
each pair of states are approximately balanced, in the
sense of \eqref{eq:pairwise_balance} below. Under the
additional regularity assumption stated below, the term
$n^{-1/3}$ improves to $n^{-1/2}$.

Importantly, for two-state systems ($K=2$), the balance
condition \eqref{eq:pairwise_balance} below holds automatically. Indeed, for
$r\in\mathcal P(\mathbb S)$,
$$
 \bnF(r)
 =\frac{r_1\nH_{1,2}(r)-r_2\nH_{2,1}(r)}{\alpha(n)}
   \bar{\canvec}_{1,2},
$$
and hence
$$
 |r_1\nH_{1,2}(r)-r_2\nH_{2,1}(r)|
 =\frac{\alpha(n)}{\sqrt2}\|\bnF(r)\|
 \leq\frac{\|\bar F\|_\infty}{\sqrt2}\alpha(n).
$$
Thus \eqref{eq:pairwise_balance} holds with
$C_{\mathrm{bal}}=\|\bar F\|_\infty/\sqrt2$.

\begin{theorem}[Cancellation of opposite jump directions]
\label{thm:critical_pairwise}
In addition to Assumption~\ref{assum:critical_regime_PRM},
suppose that the extensions $\bnF$ and $\crdiff^{(n)}$
belong to $C_b^4$ and suppose \eqref{eq:critical_uniform_C3_PRM} holds also for $j=4$.
Suppose also that there exists $C_{\mathrm{bal}}<\infty$,
independent of $n$, such that
\begin{equation}\label{eq:pairwise_balance}
 \sup_{r\in\mathcal P(\mathbb S)}
 \sum_{\substack{k,l\\l>k}}|r_k\nH_{k,l}(r)-r_l\nH_{l,k}(r)|
 \leq C_{\mathrm{bal}}\alpha(n),
 \qquad n\geq1.
\end{equation}
Then, for every $T>0$ and $n\geq1$, the first and second
bounds in \eqref{eq:critical_dK_bound_PRM} improve to
$C'_{\mathrm{SDE},0}(T)n^{-1/2}$ and
$C'_{\mathrm{SDE},1}(T)
 \bigl[n^{-1/2}+|a_n^{-1}-a_\infty^{-1}|\bigr]$,
respectively.
Here $C'_{\mathrm{SDE},0}(T)$ and
$C'_{\mathrm{SDE},1}(T)$ are finite constants depending
only on the quantities in Theorem \ref{thm:critical_dK_PRM},
$C_{\mathrm{bal}}$,
$\|D^4\bar F\|_{\mathrm{op},\infty}$, and
$\|D^4\crdiff\|_{\mathrm{op},\infty}$.
\end{theorem}

Theorem~\ref{thm:critical_pairwise} is proved in
Section~\ref{sec:critical_pairwise_proof}.

In many models, the diffusion arising in the critical
regime is of Wright--Fisher type. Its diffusion coefficient
is not differentiable at the boundary of the simplex,
so the regularity condition on $\crdiff^{(n)}$ in
Assumption \ref{assum:critical_regime_PRM} does not hold.
Below we give conditions under which quantitative
approximation bounds can still be obtained using backward
regularity estimates on the simplex.

Let $d=K-1$ and put
$$
 E\dfeq\left\{u\in[0,1]^d:\sum_{i=1}^d u_i\leq1\right\}.
$$
For $u\in E$, define the map $\mathfrak{r}: E \to \mathcal{P}(\mathbb S)$ by
$$
 \mathfrak{r}(u)\dfeq\left(u_1,\ldots,u_d,1-\sum_{i=1}^d u_i\right).
$$
Here $C^m(E)$ consists of restrictions of $C^m$ functions
on a neighborhood of $E$. For functions defined on $E$,
all derivative norms are taken over $E$.

Recall that, by the definition of $\bnF$,
\begin{equation}\label{eq:WF_drift_boundary_PRM}
 \sum_{i=1}^K\bnF_i(r)=0,
 \qquad
 \bnF_i(r)\geq0\quad\text{if }r_i=0,
 \qquad r\in\mathcal P(\mathbb S),\quad i=1,\ldots,K.
\end{equation}
For the Wright--Fisher approximation, we assume that
\begin{equation}\label{eq:WF_drift_regularity_PRM}
 \sup_{n\geq1}
 \|\bnF\circ \mathfrak{r}\|_{C^6(E;\R^K)}<\infty.
\end{equation}
For $r\in\mathcal P(\mathbb S)$, define
$$
 \crdiff_{\mathrm{WF}}(r)
 \dfeq
 \left[2\bigl(\operatorname{diag}(r)-rr^{\mathsf T}\bigr)
 \right]^{1/2},
$$
where the symmetric nonnegative definite square root is
used.

For $c>0$ and $z\in\mathcal P(\mathbb S)$, let
$\crlim_{\mathrm{WF},z}^{(n,c)}$ denote the solution of
\begin{align}\label{eq:WF_diffusion_SDE_PRM}
\begin{aligned}
 d\crlim_{\mathrm{WF},z}^{(n,c)}(t)
 ={}&
 \bnF\bigl(\crlim_{\mathrm{WF},z}^{(n,c)}(t)\bigr)\,dt+\frac{1}{c}
 \crdiff_{\mathrm{WF}}
 \bigl(\crlim_{\mathrm{WF},z}^{(n,c)}(t)\bigr)\,dW(t),
 \quad
 \crlim_{\mathrm{WF},z}^{(n,c)}(0)=z,
\end{aligned}
\end{align}
where $W$ is a $K$-dimensional Brownian motion.
Under \eqref{eq:WF_drift_boundary_PRM} and
\eqref{eq:WF_drift_regularity_PRM}, this SDE has a
$\mathcal P(\mathbb S)$-valued weak solution that is unique
in law; see \cite[Chapter 8, Theorem 2.8]{EthierKurtz1986}.
For notational convenience, write
$$
 \crlim_{\mathrm{WF}}^{(n,c)}
 \equiv
 \crlim_{\mathrm{WF},\mu_0^{(n)}}^{(n,c)}.
$$

The generator of the first $d$ coordinates is given by
\begin{equation}\label{eq:WF_generator_PRM}
 (\mathcal B_{\mathrm{WF},c}^{(n)}h)(u)
 \dfeq
 \sum_{i=1}^d\bnF_i(\mathfrak{r}(u))\partial_i h(u)
 +\frac1{c^2}\sum_{i,j=1}^d
 u_i(\delta_{ij}-u_j)\partial_{ij}h(u),
 \qquad h\in C^2(E;\R).
\end{equation}

\begin{theorem}[Wright--Fisher diffusion approximation]
\label{thm:critical_WF_dK_PRM}
Suppose that Assumption \ref{assum:scaling-1}, Assumption \ref{assum:critical_regime_PRM}:\ref{item:an-a_infty} and
\eqref{eq:WF_drift_regularity_PRM} hold.
Put
\begin{equation}\label{eq:WF_covariance_error_PRM}
 \varepsilon_n\dfeq
 \sup_{r\in\mathcal P(\mathbb S)}
 \left\|
 \bnG(r)-2\bigl(\operatorname{diag}(r)-rr^{\mathsf T}\bigr)
 \right\|_{\mathrm{tr}}.
\end{equation}
Then, for every $T>0$, there exist finite constants
$C_{\mathrm{WF},0}(T)$ and $C_{\mathrm{WF},1}(T)$,
depending only on $a_*,T,K$,
$\|\bar G\|_{\mathrm{tr},\infty}$, and the uniform bounds
on $D^j(\bnF\circ \mathfrak{r})$, $j=1,2,3$, such that, for every
$n\geq1$,
\begin{align}\label{eq:critical_WF_dK_bound_PRM}
\begin{aligned}
 &\sup_{t\leq T}\wass_1\bigl(
 \operatorname{Law}(\bnmu_t),
 \operatorname{Law}(\crlim_{\mathrm{WF}}^{(n,a_n)}(t))
 \bigr) \leq\ 
 C_{\mathrm{WF},0}(T)
 \bigl[n^{-1/3}+\varepsilon_n^{1/2}\bigr],\\
 &\sup_{t\leq T}\wass_1\bigl(
 \operatorname{Law}(\bnmu_t),
 \operatorname{Law}(\crlim_{\mathrm{WF}}^{(n,a_\infty)}(t))
 \bigr) \leq \ 
 C_{\mathrm{WF},1}(T)
 \bigl[n^{-1/3}+\varepsilon_n^{1/2}
       +|a_n^{-2}-a_\infty^{-2}|^{1/2}\bigr].
\end{aligned}
\end{align}

If, in addition, \eqref{eq:pairwise_balance} holds, then
$n^{-1/3}$ can be replaced by $n^{-1/2}$ in both bounds,
with constants $C'_{\mathrm{WF},0}(T)$ and
$C'_{\mathrm{WF},1}(T)$ that also depend on
$C_{\mathrm{bal}}$ and the uniform bound on
$D^4(\bnF\circ \mathfrak r)$. In particular, this improvement
holds whenever $K=2$.
\end{theorem}
Theorem~\ref{thm:critical_WF_dK_PRM} is proved in Section \ref{sec:WF-diff-proof}.

\begin{remark}\label{rem:random_initial_PRM}
The above results are stated for deterministic
$\mu_0^{(n)}$. They remain valid, with the same constants,
when $\mu_0^{(n)}$ is random, provided that the initial
configuration is independent of the driving PRMs $\xi_{i,k,l}$ and the comparison Brownian motions $W$.

Indeed, all estimates are uniform in the initial empirical
measure. Conditioning on $\mu_0^{(n)}$ and then taking
expectations gives the moment and concentration bounds.
For the diffusion approximations, note that, for any
integrable $\R^K$-valued random vectors $U,V$,
$$
 \wass_1\bigl(\operatorname{Law}(U),\operatorname{Law}(V)\bigr)
 \leq
 \EE\left[
 \wass_1\bigl(
 \operatorname{Law}(U\mid\mu_0^{(n)}),
 \operatorname{Law}(V\mid\mu_0^{(n)})
 \bigr)
 \right].
$$
The conditional Wasserstein bounds therefore give the
same unconditional bounds.

Note that in this case $\rho^{(n)}$, defined by \eqref{eq:ODE-approx}
with $\rho_0^{(n)}=\mu_0^{(n)}$, is a solution of the random ODE, and $Y_{0,0}^{(n)}$ is Gaussian conditionally on $\mu_0^{(n)}$.
\end{remark}

\section[ODE regime: approximation, fluctuations and concentration]
{The ODE regime: Proofs}
\label{sec:ODE_regime_PRM}
 We now prove the results in the ODE regime.

\subsection{ODE approximation}
\label{sec:ode_moment_proof}

\begin{proof}[Proof of Theorem~\ref{th:ODE-approx}]

\ref{item:limit-is-prob}
First note that for $r \in \mathcal P(\mathbb S)$
\[
\begin{aligned}
& \sum_{i=1}^K \bnF_i(r)
=\mathbf 1^\top\bnF(r)=0,\\
& r_i=0
\ \implies\
\bnF_i(r)
=\canvec_i^\top\bnF(r)
=\alpha(n)^{-1}
  \sum_{\substack{k=1\\k\neq i}}^K r_k\nH_{k,i}(r)
\geq0.
\end{aligned}
\tag{*}
\]
For a step size $h$, define the one step Euler map $\phi_h: \mathbb U \to \R^K$ by 
$$\phi_h(u) = u+ \bnF(u)h.$$
We first establish the following claim.

{\em Claim:} For sufficiently small $h$, $\phi_h \Big|_{\mathcal P(\mathbb S)}: \mathcal P(\mathbb S) \to \mathcal P(\mathbb S).$

To see this, let $r \in \mathcal P(\mathbb S)$, and for coordinates $i \neq j$, define 
$$r^{(i,j)} = r+ r_i \bar{\canvec}_{i,j} \equiv r+r_i(\canvec_j - \canvec_i)$$
Clearly, $r^{(i,j)} \in  \mathcal P(\mathbb S)$, and $\|r - r^{(i,j)}\| = \sqrt 2 r_i$. Consequently, by $(*)$ and Lipschitz continuity of $\bnF$, we have for any coordinate $i =1,2\hdots, K$
\[
\bnF_i(r)
\geq \bnF_i(r^{(i,j)})-\|D\bar F\|_{\mathrm{op},\infty}\|r-r^{(i,j)}\|
\geq -\sqrt2\,\|D\bar F\|_{\mathrm{op},\infty}r_i.
\]
It follows for any $h >0$ with $\sqrt2 \|D\bar F\|_{\mathrm{op},\infty}h\leq 1$ 
\begin{align*}
  (\phi_h(r))_i =&\  r_i + \bnF_i(r)h \geq (1-\sqrt 2 \|D\bar F\|_{\mathrm{op},\infty} h)r_i \geq 0,\\ \quad \sum_{i=1}^K(\phi_h(r))_i =&\ \sum_{i=1}^K r_i + h \sum_{i=1}^K \bnF_i(r)= 1,  
\end{align*}
which proves the claim.
Denote the $m$-composition of $\phi_h$ as $\phi^{(\circ m)}_h$.
By the convergence of the Euler scheme, $\sup_{t\leq T}\|\rho_t^{(n)} -\phi^{(\circ \lfloor t/h\rfloor)}_h(\mu^{(n)}_0)\| \to 0 $ as $h \downarrow 0$. Since for $h >0$ with $\sqrt2 \|D\bar F\|_{\mathrm{op},\infty}h\leq 1$,  each $\phi^{(\circ \lfloor t/h\rfloor)}_h(\mu^{(n)}_0) \in \mathcal P(\mathbb S)$ by the claim, and  $\mathcal P(\mathbb S) \subset \R^K$ is closed, we conclude that $\rho_t^{(n)} \in \mathcal P(\mathbb S)$ for $t \in [0,T]$. \\

\noindent
\ref{item:ODE-approx} Notice that \eqref{eq:ODE-approx},
\eqref{eq:scaled_characteristics_main} and the Lipschitz property of $\bnF$ give
\[
\|\bnmu_t-\rho_t^{(n)}\|
 \leq \ \|D\bar F\|_{\mathrm{op},\infty}\int_0^t\|\bnmu_s
                  -\rho_s^{(n)}\|\,ds
   +\|\bnmart(t)\|
\]
It follows from Gronwall's inequality that
\begin{equation}\label{eq:gronwall-pathwise}
 \sup_{t\leq T}\|\bnmu_t-\rho_t^{(n)}\|
 \leq e^{\|D\bar F\|_{\mathrm{op},\infty}T}\sup_{t\leq T}\|\bnmart(t)\|.
\end{equation}
Now since $\Delta \bnmart(s) = \bar{\canvec}_{k,l}/n$, Lemma \ref{lem:burkholder_rosenthal_PRM}, Assumption \ref{assum:scaling-1}, and  \eqref{eq:scaled_empirical_covariation}  imply that
\[
\begin{aligned}
\EE\!\left[
 \sup_{s\le T}\|\bnmart(s)\|^p\right]
&\leq C_{\mathrm{BR}}(p)\left[
 \left(\frac{T\|\bar G\|_{\mathrm{tr},\infty}}{a_n^2}\right)^{p/2}
 +\left(\frac{\sqrt2}{n}\right)^p
 \right]\\
&\leq C_{\mathrm{BR}}(p)(T\|\bar G\|_{\mathrm{tr},\infty} \vee 2)^p\bigl(a_n^{-p}+n^{-p}\bigr).
\end{aligned}
\]
The assertion now follows from \eqref{eq:gronwall-pathwise}.

\end{proof}

\subsection{Exponential concentration}
\label{sec:concentration_proof}

The following lemma is key to the proof of Theorem \ref{thm:concentration_PRM}.

\begin{lemma}[Martingale concentration]
\label{lem:vector_martingale_conc_PRM}
For every $n\geq1$ and $T,u>0$,
\begin{equation}\label{eq:vector_martingale_tail_conc_PRM}
 \PP\left(
 \sup_{t\leq T}\|\bnmart(t)\|\geq u
 \right)
 \leq2\exp\left\{
 -\frac{u^2a_n^2}
 {2\bigl(\|\bar G\|_{\mathrm{tr},\infty}T+\sqrt2\,u a_n^2/(3n)\bigr)}
 \right\}.
\end{equation}
\end{lemma}

\begin{proof}
Fix $n\geq1$ and $T>0$. Note that \eqref{eq:scaled_characteristics_main} and \eqref{eq:scaled-drft-quad-bdd} implies that 
\begin{align}\label{eq:det-mart-bd}
\sup_{t\leq T}\|\bnmart(t)\|
\leq\sqrt2+T\|\bar{F}\|_{\infty}.
\end{align}

If $\|\bar G\|_{\mathrm{tr},\infty}=0$, then $\bnmart\equiv0$
and the assertion is immediate. Assume henceforth that
$\|\bar G\|_{\mathrm{tr},\infty}>0$. For $0<\lambda<3n/\sqrt2$, set
$$
\psi_\lambda(x)=\cosh(\lambda\|x\|),\qquad x\in\R^K.
$$

Clearly, $\psi_\lambda\in C^\infty(\R^K)$. We next need estimates on its derivatives. Note that $D\psi_\lambda(0)=0$, and for $x\neq0$,
$$
D\psi_\lambda(x)[z] = \nabla^\top \psi_\lambda(x)z
=\lambda\sinh(\lambda\|x\|)
  \frac{\langle x,z\rangle}{\|x\|},
$$
Consequently,

for any $R>0$,
\begin{align}\label{eq:derivative-psi-bd}
\sup_{\|x\|\leq R}\|D\psi_\lambda(x)\|
= \sup_{\|x\|\leq R}\lambda\sinh(\lambda\|x\|)  = \lambda\sinh(\lambda R)<\infty.
\end{align}
Next note that the second derivative operator is given by
\begin{align*}
D^2\psi_\lambda(x)[z,z] \equiv z^\top D^2\psi_\lambda(x)z 
={}&\lambda^2\cosh(\lambda\|x\|)
       \frac{\langle x,z\rangle^2}{\|x\|^2}\\
&+\frac{\lambda\sinh(\lambda\|x\|)}{\|x\|}
  \left(\|z\|^2-\frac{\langle x,z\rangle^2}{\|x\|^2}\right), \quad x \neq 0.
\end{align*}
Since $\sinh y\leq y\cosh y$ for $y\geq0$, and
$D^2\psi_\lambda(0)[z,z]=\lambda^2\|z\|^2$, we obtain
\[
0\leq D^2\psi_\lambda(x)[z,z]
\leq\lambda^2\psi_\lambda(x)\|z\|^2,
\qquad x,z\in\R^K.
\]
Moreover,
$\psi_\lambda(x+\theta z)
\leq e^{\lambda\theta\|z\|}\psi_\lambda(x)$
for $0\leq\theta\leq1$. Thus, by Taylor's expansion, for $\|z\|\leq\sqrt2/n$,
\begin{align} \label{eq:cosh_remainder_conc_PRM}
\begin{aligned}
0\leq{}&
\psi_\lambda(x+z)-\psi_\lambda(x)-D\psi_\lambda(x)[z]
\ = \int_0^1(1-\theta)
        D^2\psi_\lambda(x+\theta z)[z,z]\,d\theta \\
\leq&\ \lambda^2\psi_\lambda(x)\|z\|^2
        \int_0^1(1-\theta)e^{\lambda\theta\|z\|}\,d\theta
\ = \ \psi_\lambda(x)
       \bigl(e^{\lambda\|z\|}-1-\lambda\|z\|\bigr)\\
\leq&\ 
\frac{\lambda^2}{2(1-\lambda\sqrt2/(3n))}
        \psi_\lambda(x)\|z\|^2.
\end{aligned}
\end{align}
Here the last inequality follows from
$e^y-1-y\leq y^2/\big(2(1-y/3)\big)$ for $0\leq y<3$.

Next, applying  It\^o's lemma to \eqref{eq:scaled_empirical_martingale_PRM} we get
\begin{align} \label{eq:Ito_cosh_bmart_conc_PRM}
\begin{aligned}
\psi_\lambda(\bnmart(t))
=&\ 1
+\sum_{i=1}^n\sum_{\substack{k,l=1\\k\neq l}}^K
  \int_{[0,\infty)\times(0,t]}
  \left[
  \psi_\lambda\left(\bnmart(s-)+\frac{\bar{\canvec}_{k,l}}n\right)
  -\psi_\lambda(\bnmart(s-))
  \right]
\\
&\qquad\times
  \indic_{\{\bnX_i(s-)=k\}}
  \indic_{[0,\nH_{k,l}(\bnmu_{s-})]}(w)
  \,\widetilde\xi_{i,k,l}
       \bigl(dw\times\alpha(n)^{-1}ds\bigr)
\\
&\ +\frac{n}{\alpha(n)}\int_0^t
  \sum_{\substack{k,l=1\\k\neq l}}^K\bnmu_{s-}(k)\nH_{k,l}(\bnmu_{s-})\Bigg[
  \psi_\lambda\left(\bnmart(s-)+\frac{\bar{\canvec}_{k,l}}n\right)
  \\
&\qquad -\psi_\lambda(\bnmart(s-)) -\frac1nD\psi_\lambda(\bnmart(s-))[\bar{\canvec}_{k,l}]
  \Bigg]ds.
\end{aligned}
\end{align}
Denote the stochastic integral in \eqref{eq:Ito_cosh_bmart_conc_PRM}
by $\mathcal N_t^{(n)}$. Then the mean value
theorem gives
\[
\left|
\psi_\lambda\left(\bnmart(s-)+\frac{\bar{\canvec}_{k,l}}n\right)
-\psi_\lambda(\bnmart(s-))
\right|
\leq C_{\lambda,T}\frac{\|\bar{\canvec}_{k,l}\|}{n},
\qquad 0<s\leq T,
\]
where by \eqref{eq:det-mart-bd} and \eqref{eq:derivative-psi-bd} the deterministic constant $C_{\lambda,T} $ can be taken as  $$C_{\lambda,T} \equiv \sup \lf\{\|D\psi_\lambda(x)\|: \|x\|<2\sqrt2+T\|\bar{F}\|_{\infty}\ri\} = \lambda\sinh(\lambda(2\sqrt2+T\|\bar{F}\|_{\infty})) < \infty.$$ 
Also note that
\begin{align}\label{eq:pred-cov}
\sum_{\substack{k,l=1\\k\neq l}}^K
\bnmu_{s}(k)\nH_{k,l}(\bnmu_{s})
\|\bar{\canvec}_{k,l}\|^2
=\beta(n)\trace\bigl(\bnG(\bnmu_{s})\bigr)
\ \leq\ \beta(n)\|\bar G\|_{\mathrm{tr},\infty}.
\end{align}
The predictable quadratic variation of $\mathcal N^{(n)}$
therefore satisfies
\begin{align*}
\EE\left[\langle\mathcal N^{(n)}\rangle_T\right]
={}&\frac{n}{\alpha(n)}\EE\int_0^T
\sum_{\substack{k,l=1\\k\neq l}}^K\bnmu_{s-}(k)\nH_{k,l}(\bnmu_{s-})
\\
&\quad\times\left|
\psi_\lambda\left(\bnmart(s-)+\frac{\bar{\canvec}_{k,l}}n\right)
-\psi_\lambda(\bnmart(s-))
\right|^2\,ds
\\
\leq&\ \frac{C_{\lambda,T}^2}{n\alpha(n)}
\EE\int_0^T\sum_{\substack{k,l=1\\k\neq l}}^K
\bnmu_{s}(k)\nH_{k,l}(\bnmu_{s})
\|\bar{\canvec}_{k,l}\|^2\,ds
\leq\ \frac{C_{\lambda,T}^2
\|\bar G\|_{\mathrm{tr},\infty}T}{a_n^2}<\infty.
\end{align*}
Consequently, $\mathcal N^{(n)}$ is a square-integrable martingale
on $[0,T]$. 

Hence \eqref{eq:Ito_cosh_bmart_conc_PRM}, \eqref{eq:pred-cov} and \eqref{eq:cosh_remainder_conc_PRM} imply that
\begin{align*}
&\quad \EE\bigl[\psi_\lambda(\bnmart(t))\bigr]\\
&\leq 1
+\frac{\lambda^2n}
       {2\alpha(n)(1-\lambda\sqrt2/(3n))}
  \int_0^t
  \EE\Biggl[
    \psi_\lambda(\bnmart(s-))
    \sum_{\substack{k,l=1\\k\neq l}}^K
    \bnmu_{s-}(k)\nH_{k,l}(\bnmu_{s-})
    \|\bar{\canvec}_{k,l}\|^2 n^{-2}
  \Biggr]\,ds\\
&\leq 1
+\frac{\lambda^2\|\bar G\|_{\mathrm{tr},\infty}}
       {2a_n^2(1-\lambda\sqrt2/(3n))}
  \int_0^t
  \EE\bigl[\psi_\lambda(\bnmart(s))\bigr]\,ds.
\end{align*}
Gronwall's inequality now shows that
\[
\EE\bigl[\psi_\lambda(\bnmart(T))\bigr]
\leq\exp\left\{
\frac{\lambda^2\|\bar G\|_{\mathrm{tr},\infty}T}
     {2a_n^2(1-\lambda\sqrt2/(3n))}
\right\}.
\]
Since $\psi_\lambda$ is convex,
$\psi_\lambda(\bnmart)$ is a nonnegative submartingale.
Doob's maximal inequality and $\cosh(\lambda u)\geq e^{\lambda u}/2$
therefore imply
\begin{align*}
\PP\left(\sup_{t\leq T}\|\bnmart(t)\|\geq u\right)
&\leq
\frac{\EE[\psi_\lambda(\bnmart(T))]}{\cosh(\lambda u)}\\
&\leq2\exp\left\{
-\lambda u+
\frac{\lambda^2\|\bar G\|_{\mathrm{tr},\infty}T}
     {2a_n^2(1-\lambda\sqrt2/(3n))}
\right\}.
\end{align*}
The assertion now follows upon taking
\[
\lambda=
\frac{u a_n^2}
     {\|\bar G\|_{\mathrm{tr},\infty}T
       +\sqrt2\,u a_n^2/(3n)}
\in(0,3n/\sqrt2).
\]

\end{proof}

\begin{proof}[Proof of Theorem~\ref{thm:concentration_PRM}]
We note that \eqref{eq:gronwall-pathwise} and Lemma~\ref{lem:vector_martingale_conc_PRM} give
\begin{align*}
 \PP\left(
 \sup_{t\leq T}\|\bnmu_t-\rho_t^{(n)}\|\geq\delta
 \right)
 &\leq\PP\left(
 \sup_{t\leq T}\|\bnmart(t)\|\geq\delta_T\right)\\
 &\leq2\exp\left\{
 -\frac{\delta_T^2a_n^2}
 {2\bigl(\|\bar G\|_{\mathrm{tr},\infty}T+\sqrt2\,\delta_Ta_n^2/(3n)\bigr)}
 \right\},
\end{align*}
where  $\delta_T=\delta e^{-\|D\bar F\|_{\mathrm{op},\infty}T}$.  
This proves the first inequality in \eqref{eq:concentration_main_PRM}.

For the second inequality, simply note that 
\[
 \frac{\delta_T^2a_n^2}
 {2\bigl(\|\bar G\|_{\mathrm{tr},\infty}T+\sqrt2\,\delta_Ta_n^2/(3n)\bigr)}
 \geq
 \frac{\delta_T^2}
 {2\bigl(\|\bar G\|_{\mathrm{tr},\infty}T+\sqrt2\,\delta_T/3\bigr)}
 (a_n^2\wedge n).
\]
Indeed, if $a_n^2\leq n$, we use $a_n^2/n\leq1$ in the denominator.
If $a_n^2>n$, we divide the numerator and denominator by $a_n^2/n$
and use $n/a_n^2\leq1$.
\end{proof}

\subsection{Linear diffusion approximation}\label{sec:gaussian_proof}

We need the following regularity estimates. Recall that $\mathcal{A}_s^{(n)}$ was defined in \eqref{eq:Lscript_PRM}.
Let $\Phi^{(n)}(t,s)$ be the fundamental matrix satisfying
\begin{equation}\label{eq:Phi_PRM_def}
 \partial_t\Phi^{(n)}(t,s)
 =D\bnF(\rho_t^{(n)})\Phi^{(n)}(t,s),
 \qquad \Phi^{(n)}(s,s)=I_K.
\end{equation}

\begin{lemma}[The linear backward equation]
\label{lem:backward_reg_PRM}
For $f\in C_b^4(\mathbb R^K,\R)$ and $0\leq s\leq t$, define
\begin{equation}\label{eq:backward-explicit-PRM}
 g_{s,t}(x)=\EE\left[f(Y_{s,x}^{(n)}(t))\right].
\end{equation}
This is the unique bounded classical solution, continuous up to its
terminal condition, of
\begin{equation}\label{eq:bke_PRM}
 (\partial_s+\mathcal A_s^{(n)})g_{s,t}=0,
 \qquad g_{t,t}=f.
\end{equation}
For $m=1,2,3$ and $0\leq s\leq t\leq T$,
\begin{align}\label{eq:backward_reg_PRM}
\begin{aligned}
 \|D^m g_{s,t}\|_{\mathrm{op},\infty}
 &\leq
 \|\Phi^{(n)}(t,s)\|_{\mathrm{op}}^m
 \|D^m f\|_{\mathrm{op},\infty}\\
 &\leq
 \exp\left(m\|D\bar F\|_{\mathrm{op},\infty}(t-s)\right)
 \|D^m f\|_{\mathrm{op},\infty},
\end{aligned}
\end{align}
where $\|D\bar F\|_{\mathrm{op},\infty}$ is as in \eqref{eq:barF_G_reg_PRM}.
\end{lemma}
The proof of Lemma \ref{lem:backward_reg_PRM} is in Appendix \ref{Regularity_for_ode_proof}.

The next comparison uses a test function $\psi$ on fluctuation
space and its composition $\Psi_s^{(n)}$ with the map from
configurations to fluctuations. The accelerated configuration
generator $\bngen$ acts on $\Psi_s^{(n)}$, while the diffusion
generator $\mathcal A_s^{(n)}$ acts on $\psi$.

\begin{lemma}[Generator comparison]\label{lem:gen_comp_PRM}
Let $\psi\in C_b^3(\mathbb R^K;\mathbb R)$ and, for fixed
$s\in[0,T]$, define $\Psi_s^{(n)}:\mathbb S^n\to\mathbb R$ by
$$
\Psi_s^{(n)}(x)
=\psi\bigl(a_n(\nden(x)-\rho_s^{(n)})\bigr),
$$
and set $y=a_n(\nden(x)-\rho_s^{(n)})$.
Then 
\begin{align}
\label{eq:gen_comp_bound_PRM}
 \bigl|(\bngen \Psi_s^{(n)})(x)
 &-a_n\bar F^{(n)}(\rho_s^{(n)})\cdot D\psi(y)
  -\mathcal A_s^{(n)}\psi(y)\bigr|\notag\\
 &\leq C_{\mathrm{gen}}\left(
 a_n^{-1}\|D\psi\|_{\mathrm{op},\infty}\|y\|^2
 +a_n^{-1}\|D^2\psi\|_{\mathrm{op},\infty}\|y\|
 +\frac{a_n}{n}\|D^3\psi\|_{\mathrm{op},\infty}\right)
\end{align}
where the constant $C_{\mathrm{gen}} \equiv C_{\mathrm{gen}}(\|\bar G\|_{\mathrm{tr},\infty}, L_{DF}, L_G)$ is given by
\begin{align}\label{eq:gen_comp_const}
C_{\mathrm{gen}} =\frac{L_{DF}}2
\vee\frac{L_G}2
\vee\frac{\sqrt2}{6}\|\bar G\|_{\mathrm{tr},\infty}.
\end{align}
\end{lemma}

\begin{proof}
Fix $s\in[0,T]$ and $x\in\mathbb S^n$, and put
$$
r=\nden(x)=\rho_s^{(n)}+\frac{y}{a_n}.
$$
By \eqref{eq:scaled_empirical_generator_PRM},
\begin{equation}\label{eq:nalphaG_hat_h}
 \bngen \Psi_s^{(n)}(x)
 =\frac n{\alpha(n)}
  \sum_{\substack{k,l=1\\k\neq l}}^K r_k\nH_{k,l}(r)
  \bigl[\psi\lf(y+\frac{a_n}{n}\bar{\canvec}_{k,l}\ri)-\psi(y)\bigr].
\end{equation}
Taylor's formula gives
\begin{align}\label{eq:Taylor3_PRM}
 \psi(y+\Delta)-\psi(y)
 &=D\psi(y)[\Delta]
   +\frac12D^2\psi(y)[\Delta,\Delta]
   +R_\psi(y,\Delta),
\end{align}
where
\begin{align}\label{eq:taylor-rem-0}
  |R_\psi(y,\Delta)|
 \leq\frac16\|D^3\psi\|_{\mathrm{op},\infty}\|\Delta\|^3.
\end{align}
Now notice that by the definition of $\bnF$ and $\bnG$ (see \eqref{eq:scaled-drft-quad-def-1}) the linear and quadratic terms in \eqref{eq:nalphaG_hat_h}
are, respectively,
\begin{align*}
 \frac n{\alpha(n)}
 \sum_{\substack{k,l=1\\k\neq l}}^K
 r_k\nH_{k,l}(r)D\psi(y)\lf[\frac{a_n}{n}\bar{\canvec}_{k,l}\ri]
 &=a_n\bnF(r)\cdot D\psi(y),\\
 \frac n{2\alpha(n)}
 \sum_{\substack{k,l=1\\k\neq l}}^K
 r_k\nH_{k,l}(r)
 D^2\psi(y)\lf[\frac{a_n}{n}\bar{\canvec}_{k,l},\frac{a_n}{n}\bar{\canvec}_{k,l}\ri]
 &=\frac{a_n^2}{2n\alpha(n)} \trace\lf(G^{(n)}(r)D^2\psi(y)\ri)\\
 &=\frac{1}{2} \trace\lf(\bnG(r)D^2\psi(y)\ri).
\end{align*}

Next, notice that since 
$\|\bar{\canvec}_{k,l}\|=\sqrt2$, we have by \eqref{eq:taylor-rem-0} and \eqref{eq:trace_total_rate}
\begin{align*}
R \dfeq& \ \frac n{\alpha(n)}
  \sum_{\substack{k,l=1\\k\neq l}}^K r_k\nH_{k,l}(r) \lf|R_\psi\lf(y,\frac{a_n}{n}\bar{\canvec}_{k,l}\ri)\ri|\\
  &\leq \ \frac{n}{6\alpha(n)}\|D^3\psi\|_{\mathrm{op},\infty}
   \sum_{\substack{k,l=1\\k\neq l}}^K
   r_k\nH_{k,l}(r)\lf\|\frac{a_n}{n}\bar{\canvec}_{k,l}\ri\|^3\\
  \leq &\  \frac{\sqrt 2 a_n^3}{6n^2\alpha(n)}\|D^3\psi\|_{\mathrm{op},\infty}
   \sum_{\substack{k,l=1\\k\neq l}}^K
   r_k\nH_{k,l}(r)\lf\|\bar{\canvec}_{k,l}\ri\|^2 = \frac{\sqrt2\,a_n^3\beta(n)}{6n^2\alpha(n)}
 \|D^3\psi\|_{\mathrm{op},\infty}\operatorname{tr}\bar G^{(n)}(r)\\
&\leq \frac{\sqrt2\,a_n}{6n}
 \|D^3\psi\|_{\mathrm{op},\infty}\|\bar G\|_{\mathrm{tr},\infty},
\end{align*}
where we recall that $a_n^2 \equiv n\alpha(n)/\beta(n)$.
Next observe that $r,\rho_s^{(n)}\in\mathcal P(\mathbb S)$, and 
the segment joining them lies in $\mathcal P(\mathbb S)$. Therefore, by Taylor's theorem and \eqref{eq:barF_G_reg_PRM}
\begin{align}
\label{eq:drift_integral_remainder_ode}
\begin{aligned}
 &\|a_n\bigl(\bar F^{(n)}(r)-  \bar F^{(n)}(\rho_s^{(n)})\bigr)
   - D\bar F^{(n)}(\rho_s^{(n)})y\|\\
   & \qquad \qquad = \lf\|\int_0^1 \lf( D\bnF(\rho_s^{(n)}+\frac{\theta y}{a_n})
   -D\bar F^{(n)}(\rho_s^{(n)})\ri)y\ d\theta\ri\|\\
 & \qquad \qquad   \leq\ L_{DF} \|y\|^2/(2a_n)\\
 &\|\bar G^{(n)}(r)-  \bar G^{(n)}(\rho_s^{(n)})\|_{\mathrm{tr}}\ \leq  \ L_G\|y\|/a_n.
\end{aligned}
\end{align}
It now follows from the above estimates, the definition of $\mathcal A_s^{(n)}$ (see \eqref{eq:Lscript_PRM}), and the inequality $|\trace(AB)| \leq \|A\|_{\mathrm{tr}}\|B\|_{\mathrm{op}}$  that 
\begin{align*}
 &\left|
 (\bngen \Psi_s^{(n)})(x)
 -a_n\bar F^{(n)}(\rho_s^{(n)})\cdot D\psi(y)
 -\mathcal A_s^{(n)}\psi(y)
 \right|\\
 &\quad\leq
 \left\|
 a_n\bigl(\bar F^{(n)}(r)-\bar F^{(n)}(\rho_s^{(n)})\bigr)
 -D\bar F^{(n)}(\rho_s^{(n)})y
 \right\|\|D\psi\|_{\mathrm{op},\infty}\\
 &\qquad\quad
 +\frac{1}{2}
 \|\bar G^{(n)}(r)-\bar G^{(n)}(\rho_s^{(n)})\|_{\mathrm{tr}}
 \|D^2\psi\|_{\mathrm{op},\infty}+R\\
 &\quad\leq C_{\mathrm{gen}}\left(
 \frac{\|D\psi\|_{\mathrm{op},\infty}}{a_n}\|y\|^2
 +\frac{\|D^2\psi\|_{\mathrm{op},\infty}}{a_n}\|y\|
 +\frac{a_n}{n}\|D^3\psi\|_{\mathrm{op},\infty}
 \right),
\end{align*}
where $C_{\mathrm{gen}}$ is given by \eqref{eq:gen_comp_const}.
\end{proof}

\begin{proof}[Proof of Theorem~\ref{thm:kantorovich_PRM}]
Notice that, since $\rho_0^{(n)}=\bnmu_0$, we have
$\fluct_0^{(n)}=0$, and
\begin{align}\label{eq:xi_semimart_PRM}
 \fluct_t^{(n)}
 &=a_n\int_0^t
 \bigl(\bnF(\bnmu_s)-\bnF(\rho_s^{(n)})\bigr)\,ds
 +a_n\bnmart(t).
\end{align}
By Gronwall's inequality, Doob's inequality, and
the Lipschitz condition on $\bnF$,
\begin{align}\label{eq:xi_second_moment_PRM}
 \EE\left[\sup_{s\leq T}\|\fluct_s^{(n)}\|^2\right]
 &\leq
 a_n^2 e^{2T\|D\bar F\|_{\mathrm{op},\infty}}
 \EE\left[\sup_{s\leq T}\|\bnmart(s)\|^2\right]\notag\\
 &\leq
 4e^{2T\|D\bar F\|_{\mathrm{op},\infty}}
 \EE\int_0^T\trace\bigl(\bnG(\bnmu_s)\bigr)\,ds\notag\\
 &\leq
 4T e^{2T\|D\bar F\|_{\mathrm{op},\infty}}
 \|\bar G\|_{\mathrm{tr},\infty}.
\end{align}

Fix $t\leq T$ and $\psi\in C_b^\infty(\R^K;\R)$.
Let $g_{s,t}$ be the backward function in
\eqref{eq:backward-explicit-PRM} with terminal value $ g_{t,t} \equiv \psi$.
For $0\leq s\leq t$ and $x\in\mathbb S^n$, define
$$
 \Psi_{s,t}^{(n)}(x)
 \dfeq g_{s,t}\bigl(a_n(\nden(x)-\rho_s^{(n)})\bigr).
$$
By It\^o's formula, we get from \eqref{eq:xi_semimart_PRM} and \eqref{eq:scaled_empirical_martingale_PRM} that for
$0\leq u\leq t$,
\begin{align}\label{eq:ito_backward_xi_PRM}
\begin{aligned}
 g_{u,t}(\fluct_u^{(n)})
 ={}&g_{0,t}(0)+\mathcal N_{u,t}^{(n)}
 +\int_0^u(\partial_sg_{s,t})(\fluct_s^{(n)})\,ds\\
 &+a_n\int_0^u
 \bigl(\bnF(\bnmu_s)-\bnF(\rho_s^{(n)})\bigr)
 \cdot Dg_{s,t}(\fluct_s^{(n)})\,ds\\
 &+\frac n{\alpha(n)}\int_0^u
 \sum_{k\neq l}\bnmu_s(k)\nH_{k,l}(\bnmu_s)\\
 &\qquad\times
 \bigl[
 g_{s,t}\lf(\fluct_s^{(n)}+\frac{a_n}{n}\bar{\canvec}_{k,l}\ri)
 -g_{s,t}(\fluct_s^{(n)})
 -Dg_{s,t}(\fluct_s^{(n)})\lf[\frac{a_n}{n}\bar{\canvec}_{k,l}\ri]
 \bigr]\,ds,
\end{aligned}
\end{align}
where
\begin{align*}
 \mathcal N_{u,t}^{(n)}
 ={}&\sum_{i=1}^n\sum_{\substack{k,l=1\\k\neq l}}^K
 \int_{[0,\infty)\times[0,u]}
 \bigl[
 g_{s,t}\lf(\fluct_{s-}^{(n)}+\frac{a_n}{n}\bar{\canvec}_{k,l}\ri)
 -g_{s,t}(\fluct_{s-}^{(n)})
 \bigr]\\
 &\qquad\times
 \indic_{[0,\nH_{k,l}(\bnmu_{s-})]}(z)
 \indic_{\{\bnX_i(s-)=k\}}
 \widetilde\xi_{i,k,l}
 \bigl(dz\times\alpha(n)^{-1}ds\bigr).
\end{align*}
Recall that $a_n^2=n\alpha(n)/\beta(n),$ and $\bnG = \beta(n)^{-1}G^{(n)}$. Hence \eqref{eq:trace_total_rate} and \eqref{eq:backward_reg_PRM} show that
\begin{align*}
 \left\langle\mathcal N_{\cdot,t}^{(n)}\right\rangle_u
 &\leq
 \int_0^u
 \|Dg_{s,t}\|_{\mathrm{op},\infty}^2
 \trace\bigl(\bnG(\bnmu_s)\bigr)\,ds\leq
 T e^{2T\|D\bar F\|_{\mathrm{op},\infty}}
 \|\bar G\|_{\mathrm{tr},\infty}\|D\psi\|_{\mathrm{op},\infty}^2.
\end{align*}
Thus $\mathcal N_{\cdot,t}^{(n)}$ is a square-integrable
martingale with mean zero. 
 Next, recalling
$\bnF=\alpha(n)^{-1}F^{(n)}$, we see from
\eqref{eq:cov-def-comp} and
\eqref{eq:scaled_empirical_generator_PRM} that
\begin{align*}
 &\frac n{\alpha(n)}
 \sum_{\substack{k,l=1\\k\neq l}}^K\bnmu_s(k)\nH_{k,l}(\bnmu_s)
 Dg_{s,t}(\fluct_s^{(n)})
 \lf[\frac{a_n}{n}\bar{\canvec}_{k,l}\ri]
 =
 a_n\bnF(\bnmu_s)\cdot Dg_{s,t}(\fluct_s^{(n)});\\
 &\bngen\Psi_{s,t}^{(n)}(\bnX(s))
 =\frac n{\alpha(n)}
 \sum_{\substack{k,l=1\\k\neq l}}^K
 \bnmu_s(k)\nH_{k,l}(\bnmu_s)
 \lf[
 g_{s,t}\lf(\fluct_s^{(n)}
       +\frac{a_n}{n}\bar{\canvec}_{k,l}\ri)
 -g_{s,t}(\fluct_s^{(n)})
 \ri].
\end{align*}
Therefore,
$$
\begin{aligned}
 \EE[g_{u,t}(\fluct_u^{(n)})]
 ={}&g_{0,t}(0)
 +\EE\int_0^u
 \Big((\partial_sg_{s,t})(\fluct_s^{(n)})
     +(\bngen\Psi_{s,t}^{(n)})(\bnX(s))\Big)\,ds\\
 &\qquad
 -a_n\EE\int_0^u\bnF(\rho_s^{(n)})
 \cdot Dg_{s,t}(\fluct_s^{(n)})\,ds\\
 ={}&g_{0,t}(0)
 +\EE\int_0^u
 \Big((\bngen\Psi_{s,t}^{(n)})(\bnX(s))
     -(\mathcal A_s^{(n)}g_{s,t})(\fluct_s^{(n)})\Big)\,ds\\
 &\qquad
 -a_n\EE\int_0^u\bnF(\rho_s^{(n)})
 \cdot Dg_{s,t}(\fluct_s^{(n)})\,ds,
\end{aligned}
$$
where for the last equality we used the backward PDE
\eqref{eq:bke_PRM}. Notice that it also gives, by It\^o's formula,
$$
\begin{aligned}
 \EE[g_{u,t}(Y_{0,0}^{(n)}(u))]
 &=g_{0,t}(0)
 +\EE\int_0^u
 \bigl[(\partial_s+\mathcal A_s^{(n)})g_{s,t}\bigr]
 (Y_{0,0}^{(n)}(s))\,ds \ =\ g_{0,t}(0).
\end{aligned}
$$
Using the terminal condition $g_{t,t}=\psi$ gives
\begin{align}\label{eq:stein_identity_PRM}
\begin{aligned}
 \EE[\psi(\fluct_t^{(n)})]
 -\EE[\psi(Y_{0,0}^{(n)}(t))]
 ={}&\int_0^t\EE\bigl[
 (\bngen\Psi_{s,t}^{(n)})(\bnX(s))\\
 &\qquad
 -a_n\bnF(\rho_s^{(n)})\cdot Dg_{s,t}(\fluct_s^{(n)})
 -(\mathcal A_s^{(n)}g_{s,t})(\fluct_s^{(n)})
 \bigr]\,ds.
\end{aligned}
\end{align}
By Lemma~\ref{lem:gen_comp_PRM} 
\begin{align}\label{eq:stein_after_dynkin_PRM}
 &\left|
 \EE[\psi(\fluct_t^{(n)})]
 -\EE[\psi(Y_{0,0}^{(n)}(t))]
 \right|\notag\\
 &\quad\leq C_{\mathrm{gen}}\int_0^t
 \biggl(
 \frac{\|Dg_{s,t}\|_{\mathrm{op},\infty}}{a_n}
 \EE\|\fluct_s^{(n)}\|^2
 +\frac{\|D^2g_{s,t}\|_{\mathrm{op},\infty}}{a_n}
 \EE\|\fluct_s^{(n)}\|
 +\frac{a_n}{n}\|D^3g_{s,t}\|_{\mathrm{op},\infty}
 \biggr)\,ds,
 \end{align}
where  $C_{\mathrm{gen}}$ is as in \eqref{eq:gen_comp_const}.
It now follows by Lemma \ref{lem:backward_reg_PRM} and \eqref{eq:xi_second_moment_PRM} 
that 
\begin{align}\label{eq:test_fn_bound_PRM}
 \left|
 \EE[\psi(\fluct_t^{(n)})]
 -\EE[\psi(Y_{0,0}^{(n)}(t))]
 \right|
 \leq C_{\mathrm{test}}(T)\left(
 \frac{\|D\psi\|_{\mathrm{op},\infty}
       +\|D^2\psi\|_{\mathrm{op},\infty}}{a_n}
 +\frac{a_n}{n}\|D^3\psi\|_{\mathrm{op},\infty}
 \right),
\end{align}
where the constant 
$
 C_{\mathrm{test}}(T) \equiv\bar C_{\mathrm{test}}\bigl(T,
 \|D\bar F\|_{\mathrm{op},\infty}, \|\bar G\|_{\mathrm{tr},\infty},
 L_{DF},L_G \bigr)
$ can be taken as
$$
  C_{\mathrm{test}}(T)\dfeq
 1\vee\left[
 C_{\mathrm{gen}}T e^{3T\|D\bar F\|_{\mathrm{op},\infty}}
 \bigl(1\vee4T\|\bar G\|_{\mathrm{tr},\infty}\bigr)
 \right].
$$
Since both laws have finite first moments, we get by 
Lemma~\ref{lem:regularization_dK} with $\ell=3$ that \eqref{eq:dK_bound_PRM} holds with 
$ C_{\mathrm{CLT}}(T) \dfeq C_{\mathrm{reg}}(K,3)\max_{1\leq j \leq 3}C^{1/j}_{\mathrm{test}}(T).$
\end{proof}

\section{Critical SDE regime: nonlinear diffusion approximation}
\label{sec:critical_nonlinear_diffusion_PRM}
Recall that $\crlim_z^{(n,c)}$ is the solution of the SDE \eqref{eq:critical_diffusion_SDE_PRM}.

\begin{lemma}[Backward regularity]
\label{lem:critical_backward_reg_PRM}
For $\psi\in C_b^3(\R^K;\R)$, $c \in (0,\infty)$ define
$$
 \bar g_{s,t}(z)
 \dfeq\EE\left[\psi(\crlim_z^{(n,c)}(t-s))\right],
 \qquad 0\leq s\leq t.
$$
Fix $T>0$, $m \in \{1,2,3\}$ and suppose that Assumption~\ref{assum:critical_regime_PRM} holds.  Then there exists a constant  
$$ C_{\mathrm{back},m}(c,T) \equiv C_{\mathrm{back},m}(c,K,T, \bigl\{\|D^j\bar F\|_{\mathrm{op},\infty},
       \|D^j\Sigma\|_{\mathrm{op},\infty}:j=1,\ldots,m\bigr\}),$$ such that for $0\leq s\leq t\leq T$,
\begin{equation}\label{eq:critical_backward_reg_PRM}
 \|\bar g_{s,t}\|_\infty\leq\|\psi\|_\infty,
 \qquad
 \|D^m\bar g_{s,t}\|_{\mathrm{op},\infty}
 \leq C_{\mathrm{back},m}(c,T)\sum_{j=1}^m
       \|D^j\psi\|_{\mathrm{op},\infty}.
\end{equation}
Moreover, $\bar g$ satisfies the classical backward equation
\begin{equation}\label{eq:critical_bke_PRM}
 (\partial_s+\mathcal B_c^{(n)})\bar g_{s,t}=0,
 \qquad 0\leq s<t,\qquad \bar g_{t,t}=\psi.
\end{equation}
If \eqref{eq:critical_uniform_C3_PRM} holds for $j=4$ and
$\psi\in C_b^4(\R^K;\R)$, the derivative estimate also holds
for $m=4$.

The same conclusions hold for the corresponding transition function
and generator when $c$ in
\eqref{eq:critical_diffusion_SDE_PRM} is replaced by $a_n$,
with constants independent of $n$.
\end{lemma}

The proof of Lemma~\ref{lem:critical_backward_reg_PRM} is given in Appendix~\ref{Regularity_for_Sde_proof}.

For Wright--Fisher diffusions, Assumption~\ref{assum:critical_regime_PRM} fails at the boundary. For this reason, we establish backward regularity for these diffusions in the following lemma separately.

\begin{lemma}[Backward regularity for Wright--Fisher diffusions]
\label{lem:critical_WF_backward_regularity_PRM}
Suppose that \eqref{eq:WF_drift_boundary_PRM} and
\eqref{eq:WF_drift_regularity_PRM} hold.
Let $d=K-1$, and let $\mathfrak p:\R^K \to \R^d$ be the canonical projection; that is, for $z \in \R^K$,
$\mathfrak p(z)\dfeq(z_1,\ldots,z_d)$.
For $\psi\in C^4(E;\R)$ and $c>0$, define
$$
 \bar g_{s,t}(x)
 \dfeq
 \EE\left[
 \psi\left(
 \mathfrak p\bigl(
 \crlim_{\mathrm{WF},\mathfrak r(x)}^{(n,c)}(t-s)
 \bigr)
 \right)
 \right],
 \qquad 0\leq s\leq t,\quad x\in E.
$$
Then $\bar g_{\cdot,t}\in C([0,t];C^4(E;\R))$ and
satisfies the classical backward equation
\begin{equation}\label{eq:WF_bke_PRM}
 (\partial_s+\mathcal B_{\mathrm{WF},c}^{(n)})
 \bar g_{s,t}=0
 \quad\text{on }E,\qquad 0\leq s<t,
 \qquad \bar g_{t,t}=\psi.
\end{equation}
Moreover, for every $T>0$ and $m\in\{1,2,3,4\}$, there
exists a finite constant
$C_{\mathrm{back},m}^{\mathrm{WF}}(c,T)$, independent of
$n$, such that, for $0\leq s\leq t\leq T$,
\begin{equation}\label{eq:WF_backward_regularity_PRM}
 \|\bar g_{s,t}\|_\infty\leq\|\psi\|_\infty,
 \qquad
 \|D^m\bar g_{s,t}\|_{\mathrm{op},\infty}
 \leq C_{\mathrm{back},m}^{\mathrm{WF}}(c,T)
       \sum_{j=1}^m\|D^j\psi\|_{\mathrm{op},\infty}.
\end{equation}
The constants depend only on $c,K,m,T$ and the uniform
bounds on $D^j(\bnF\circ\mathfrak r)$, $j=1,\ldots,m$,
and can be chosen nonincreasing in $c$.

If $K=2$ and
$$
 \bnF_1(x,1-x)=A_n(1-x)-B_nx,
 \qquad A_n,B_n>0,
$$
then, for every $0\leq s\leq t$ and $m\in\{1,2,3,4\}$,
\begin{equation}\label{eq:WF_neutral_derivative_decay_PRM}
 \|\partial_x^m\bar g_{s,t}\|_\infty
 \leq
 e^{-m[A_n+B_n+c^{-2}(m-1)](t-s)}
 \|\psi^{(m)}\|_\infty.
\end{equation}
\end{lemma}

The proof of Lemma \ref{lem:critical_WF_backward_regularity_PRM} is in Appendix~\ref{Regularity_for_Sde_proof}.

\subsection{The empirical generator and its Taylor expansion}
\label{subsec:critical_generator_proof}

The following result is the key to proving Theorem \ref{thm:critical_dK_PRM}.

\begin{lemma}[Generator comparison]
\label{lem:critical_gen_comp_PRM}
Let $h\in C_b^3(\mathbb R^K,\R)$ and $c \in (0,\infty)$. 
Then for $x \in \mathbb S^n$, setting $r \equiv r_x = \nden(x) \in \mathcal{P}(\mathbb S)$, we have 
\begin{align}
\label{eq:critical_gen_comp_PRM}
\begin{aligned}
|\bngen(h\circ\nden)(x)-\mathcal{B}_{a_n}^{(n)}h(r)|
\leq&\ \f{\sqrt 2 \|\bar G\|_{\mathrm{tr},\infty}}{6}\f{\|D^3h\|_{\mathrm{op},\infty}}{a_n^2n}\\
|\bngen(h\circ\nden)(x)-\mathcal{B}_{c}^{(n)}h(r)|
\leq&\ \|\bar G\|_{\mathrm{tr},\infty}\biggl(\f{\|D^2h\|_{\mathrm{op},\infty}}{2}\left|a_n^{-2}-c^{-2}\right|
+\frac{\sqrt 2\|D^3h\|_{\mathrm{op},\infty}}{6a_n^2n}\biggr).
\end{aligned}
\end{align}
\end{lemma}

\begin{proof}
First note that by Taylor's expansion
\begin{align*}
h\left(r+\frac{\bar{\canvec}_{k,l}}n\right)-h(r)
=\frac{1}{n}Dh(r)[\bar{\canvec}_{k,l}]
+\frac{1}{2n^2}D^2h(r)[\bar{\canvec}_{k,l},\bar{\canvec}_{k,l}]
+\f{1}{6n^3}R_{k,l}^{(n,3)}(r),
\end{align*}
where
$$
|R_{k,l}^{(n,3)}(r)|
\leq \|D^3h\|_{\mathrm{op},\infty}
       \|\bar{\canvec}_{k,l}\|^3
=2\sqrt2\,\|D^3h\|_{\mathrm{op},\infty}.
$$
Consequently, by the definition of $\bnF$ and $\bnG$ on $\mathcal{P}(\mathbb S)$ (see \eqref{eq:scaled-drft-quad-def-1}), and using $a_n^2 =  n\alpha(n)/\beta(n)$ we have from \eqref{eq:scaled_empirical_generator_PRM}
that 
\begin{align*}
    \bngen(h\circ\nden)(x) =&\ \bar F^{(n)}(r)\cdot Dh(r) + \frac{1}{2a_n^2}\trace(\bar G^{(n)}(r)D^2h(r))+ \f{1}{6n^2\alpha(n)}\bar{R}^{(n,3)}(r),\\
    \equiv &\ \mathcal{B}_{a_n}^{(n)}h(r)+ \f{1}{6n^2\alpha(n)}\bar{R}^{(n,3)}(r)
\end{align*}
where 
$\bar{R}^{(n,3)}(r) \dfeq \sum\limits_{ \substack{k, l \\ k \neq l}}
r_k\nH_{k,l}(r)R_{k,l}^{(n,3)}(r)$ can be estimated as
\begin{align*}
  |\bar{R}^{(n,3)}(r)| \leq 2\sqrt2\|D^3h\|_{\mathrm{op},\infty}\sum_{ \substack{k, l=1 \\ k \neq l}}^Kr_k\nH_{k,l}(r) \leq \sqrt2\beta(n) \|D^3h\|_{\mathrm{op},\infty} \|\bar G\|_{\mathrm{tr},\infty}
\end{align*}
This proves the first inequality in \eqref{eq:critical_gen_comp_PRM}. For the second, simply note that the inequality $|\trace(AB)|
\leq \|A\|_{\mathrm{tr}}\|B\|_{\mathrm{op}}
$ gives
\begin{align*}
    |\mathcal{B}_{a_n}^{(n)}h(r) - \mathcal{B}_{c}^{(n)}h(r)| \leq&\ \f{|\trace(\bar G^{(n)}(r)D^2h(r))|}{2}|a_n^{-2}-c^{-2}|\\
    \leq &\ \f{\|\bar G\|_{\mathrm{tr},\infty}\|D^2h\|_{\mathrm{op},\infty}}{2}|a_n^{-2}-c^{-2}|
\end{align*}
\end{proof}

\subsection{Proof of the nonlinear diffusion approximation}
\label{sec:critical_comparison_proof}
\begin{proof}[Proof of Theorem \ref{thm:critical_dK_PRM}]
For notational convenience, denote $\crlim^{(n,c)} \equiv \crlim^{(n,c)}_{\mu^{(n)}_0}$, the solution to \eqref{eq:critical_diffusion_SDE_PRM} with initial condition $\mu^{(n)}_0$ and driven by the Brownian motion $W$ independent of $\mu^{(n)}_0$. Fix $t\leq T$ and $\psi\in C_b^\infty(\R^K;\R)$, and define 
 the backward function $\bar g_{s,t}$ by 
 $
 \bar g_{s,t}(x)
 \dfeq\EE\left[\psi(\crlim^{(n,a_n)}_x(t-s))\right],\ 0\leq s\leq t.
$

By It\^o's formula  and \eqref{eq:scaled_empirical_generator_PRM}, we have for $0\leq r\leq t$,
\begin{align*}
 \bar g_{r,t}(\bnmu_r)
 ={}&\bar g_{0,t}(\mu_0^{(n)})+\mathcal N_{r,t}^{(n)}
 +\int_0^r\partial_s\bar g_{s,t}(\bnmu_s)\,ds\\
 &+\frac n{\alpha(n)}\int_0^r
 \sum_{k\neq l}\bnmu_s(k)\nH_{k,l}(\bnmu_s)
 \left[
 \bar g_{s,t}\left(\bnmu_s+\frac{\bar{\canvec}_{k,l}}n\right)
 -\bar g_{s,t}(\bnmu_s)
 \right]\,ds,\\
 ={}&\bar g_{0,t}(\mu_0^{(n)})+\mathcal N_{r,t}^{(n)}
 +\int_0^r\lf(\partial_s\bar g_{s,t}(\bnmu_s) +  \bngen ( \bar g_{s,t} \circ \nden)(\bnX(s))\ri)\,ds
\end{align*}
where (as in the proof of Theorem~\ref{thm:kantorovich_PRM})
\begin{align*}
 \mathcal N_{r,t}^{(n)}
 \dfeq{}&\sum_{i=1}^n\sum_{k\neq l}
 \int_{[0,\infty)\times[0,r]}
 \left[
 \bar g_{s,t}\left(\bnmu_{s-}+\frac{\bar{\canvec}_{k,l}}n\right)
 -\bar g_{s,t}(\bnmu_{s-})
 \right]\\
 &\qquad\times
 \indic_{[0,\nH_{k,l}(\bnmu_{s-})]}(z)
 \indic_{\{\bnX_i(s-)=k\}}
 \widetilde\xi_{i,k,l}
 \bigl(dz\times\alpha(n)^{-1}ds\bigr)
\end{align*}
is a  square-integrable martingale.

Taking expectations  and setting $r=t$ gives
\begin{align}\label{eq:critical_dynkin_expanded}
\begin{aligned}
 \EE[\psi(\bnmu_t)]
 -\EE[\bar g_{0,t}(\mu_0^{(n)})]
 =\int_0^t\EE\left[
 \partial_s\bar g_{s,t}(\bnmu_s)
 +\bngen (\bar g_{s,t} \circ \nden)(\bnX(s))
 \right]\,ds.
\end{aligned}
\end{align}
Next, since by Lemma \ref{lem:critical_backward_reg_PRM}, $ \bar g_{s,t}$ satisfies \eqref{eq:critical_bke_PRM} with $c=a_n$ 
$$
\EE[\psi(\crlim^{(n,a_n)}(t))] \equiv \EE[\bar g_{t,t}(\crlim^{(n,a_n)}(t))] =\EE[\bar g_{0,t}(\mu_0^{(n)})],
$$
and hence by Lemma \ref{lem:critical_gen_comp_PRM}
\begin{align}\label{eq:critical_Stein_identity_PRM}
\begin{aligned}
 &\Big|\EE[\psi(\bnmu_t)]-  \EE[\psi(\crlim^{(n,a_n)}(t))] \Big| \\
 &\leq \int_0^t\EE\left[\lf|
 \bngen ( \bar g_{s,t}\circ \nden)(\bnX(s))
 -(\mathcal B_{a_n}^{(n)}\bar g_{s,t})(\bnmu_s)
 \ri|\right]\,ds\\
 &\leq\ \frac{\sqrt2\,\|\bar G\|_{\mathrm{tr},\infty}}
 {6a_n^2n}
 \int_0^t\|D^3\bar g_{s,t}\|_{\mathrm{op},\infty}\,ds \leq\  \bar C_{\mathrm{test}}(T)a_n^{-2}n^{-1}
 \sum_{j=1}^3\|D^j\psi\|_{\mathrm{op},\infty},
\end{aligned}
\end{align}
where  
$$
\bar C_{\mathrm{test}}(T)
 \dfeq\frac{\sqrt 2 C_{\mathrm{back},3}(a_*,T)\|\bar G\|_{\mathrm{tr},\infty}T}{6}.
$$
In the above we used $\sup_{n\geq 1}C_{\mathrm{back},3}(a_n,T) \leq C_{\mathrm{back},3}(a_*,T)< \infty$, where $a_* \dfeq \inf_{n\geq 1} a_n \in (0,\infty)$. The first inequality in \eqref{eq:critical_dK_bound_PRM} now follows by Lemma \ref{lem:regularization_dK} with 

\begin{align}\label{eq:critical_exact_ratio_bound_PRM}
\begin{aligned}
C_{\mathrm{SDE},0}(T) \dfeq &\  C_{\mathrm{reg}}(K,3)
 \max_{1\leq j\leq3}
 \left(\frac{\bar C_{\mathrm{test}}(T)}{a_*^2}\right)^{1/j}.
\end{aligned}
\end{align}
For the second inequality, we work with $\crlim^{(n,a_n)}$ and $\crlim^{(n,a_\infty)}$ that are driven by the same Brownian
motion $W$.
Then note that $U^{(n)}(t) \dfeq \crlim^{(n,a_n)}(t)-\crlim^{(n,a_\infty)}(t)$ satisfies 
$$
U^{(n)}(t)=\int_0^t\left[
 \bnF(\crlim^{(n,a_n)}(s))
 -\bnF(\crlim^{(n,a_\infty)}(s))
 \right]\,ds+\int_0^t\Gamma^{(n)}(s)\,dW(s),
$$
where
\begin{align*}
 \Gamma^{(n)}(s)
 &\dfeq a_n^{-1}\crdiff^{(n)}(\crlim^{(n,a_n)}(s))
 -a_\infty^{-1}\crdiff^{(n)}(\crlim^{(n,a_\infty)}(s))\\
 &=a_n^{-1}\left[
 \crdiff^{(n)}(\crlim^{(n,a_n)}(s))
 -\crdiff^{(n)}(\crlim^{(n,a_\infty)}(s))
 \right]
 +(a_n^{-1}-a_\infty^{-1})
 \crdiff^{(n)}(\crlim^{(n,a_\infty)}(s)).
\end{align*}
Notice that 
\begin{align*}
 \trace(\Gamma^{(n)}(s)(\Gamma^{(n)}(s))^\top)
 &\leq K\|\Gamma^{(n)}(s)\|_{\mathrm{op}}^2\\
 &\leq 2K\lf(a_n^{-2}\|D\crdiff\|^2_{\mathrm{op},\infty}\|U^{(n)}(s)\|^2
 +\|\crdiff\|^2_{\mathrm{op},\infty}|a_n^{-1}-a_\infty^{-1}|^2\ri).
 \end{align*}
 It follows that 
 \begin{align*}
 &\EE\left[\sup_{s\leq t}\|U^{(n)}(s)\|^2\right]\\
 &\leq 2t\|D\bar F\|^2_{\mathrm{op},\infty}\int_0^t\EE\left[
 \left\|U^{(n)}(r)\right\|^2
 \right]\,dr + 8\int_0^t
 \EE\left[\trace(\Gamma^{(n)}(r)(\Gamma^{(n)}(r))^\top)\right]\,dr\\
 &\leq\left(2T\|D\bar F\|^2_{\mathrm{op},\infty}
 +16Ka_*^{-2}\|D\crdiff\|^2_{\mathrm{op},\infty}\right)
 \int_0^t\EE\left[\sup_{r'\leq r}\|U^{(n)}(r')\|^2\right]\,dr\\
 &\qquad+16KT \|\crdiff\|^2_{\mathrm{op},\infty}
 |a_n^{-1}-a_\infty^{-1}|^2.
\end{align*}
By Gronwall's and Cauchy-Schwarz inequality it follows 
\begin{align*}
     \EE\left[\sup_{s\leq t}\|U^{(n)}(s)\|\right] \leq C_{\mathrm{coup}}(T) |a_n^{-1}-a_\infty^{-1}|,
\end{align*}
where
$$
C_{\mathrm{coup}}(T)
\dfeq
4(KT)^{1/2}\|\crdiff\|_{\mathrm{op},\infty}
\exp\left(
T^2\|D\bar F\|_{\mathrm{op},\infty}^{\,2}
+8Ka_*^{-2}T\|D\crdiff\|_{\mathrm{op},\infty}^{\,2}
\right).
$$
 The result now follows from the triangle inequality
\begin{align*}
&\sup_{t\leq T}\wass_1\bigl(\operatorname{Law}(\bnmu_t),
 \operatorname{Law}(\crlim^{(n,a_\infty)}(t))\bigr)\\
&\leq  \sup_{t\leq T}\wass_1\bigl(\operatorname{Law}(\bnmu_t),
 \operatorname{Law}(\crlim^{(n,a_n)}(t)) \bigr)
 +\EE\left[\sup_{t\leq T}\|U^{(n)}(t)\|\right]
\end{align*}
with $C_{\mathrm{SDE},1}(T) = C_{\mathrm{SDE},0}(T) \vee C_{\mathrm{coup}}(T).$
\end{proof}

\subsection{Cancellation of opposite jump directions}
\label{sec:critical_pairwise_proof}

\begin{proof}[Proof of Theorem \ref{thm:critical_pairwise}]
Let $h\in C_b^4(\R^K;\R)$, $x\in\mathbb S^n$, and
$r=\nden(x)$. By Taylor's formula,
\begin{align}\label{eq:critical_taylor_fourth}
\begin{aligned}
 h\left(r+\frac{\bar{\canvec}_{k,l}}n\right)-h(r)
 ={}&\frac1nDh(r)[\bar{\canvec}_{k,l}]
 +\frac1{2n^2}D^2h(r)
 [\bar{\canvec}_{k,l},\bar{\canvec}_{k,l}]\\
 &+\frac1{6n^3}D^3h(r)
 [\bar{\canvec}_{k,l},\bar{\canvec}_{k,l},
  \bar{\canvec}_{k,l}]
 +R_{k,l}^{(4)}(r),
\end{aligned}
\end{align}
where, since $\|\bar{\canvec}_{k,l}\|=\sqrt2$,
$$
 |R_{k,l}^{(4)}(r)|
 \leq \frac{\|\bar{\canvec}_{k,l}\|^4}{24n^4}
       \|D^4h\|_{\mathrm{op},\infty}
 =\frac1{6n^4}\|D^4h\|_{\mathrm{op},\infty}.
$$
Substituting \eqref{eq:critical_taylor_fourth} into
\eqref{eq:scaled_empirical_generator_PRM}, the first two
terms give $(\mathcal B_{a_n}^{(n)}h)(r)$. Also, since
$\bar{\canvec}_{l,k}=-\bar{\canvec}_{k,l}$, the cubic term
satisfies
\begin{align}\label{eq:critical_signed_cubic}
\begin{aligned}
 &\frac1{6\alpha(n)n^2}
 \sum_{\substack{k,l=1\\k\neq l}}^Kr_k\nH_{k,l}(r)
 D^3h(r)[\bar{\canvec}_{k,l},\bar{\canvec}_{k,l},
         \bar{\canvec}_{k,l}]\\
 &\qquad=
 \frac1{6\alpha(n)n^2}
 \sum_{\substack{k,l\\l>k}}\bigl(r_k\nH_{k,l}(r)-r_l\nH_{l,k}(r)\bigr)
 D^3h(r)[\bar{\canvec}_{k,l},\bar{\canvec}_{k,l},
         \bar{\canvec}_{k,l}].
\end{aligned}
\end{align}
By \eqref{eq:pairwise_balance}, its absolute value is
bounded by
$$
 \frac{\sqrt2\,\|D^3h\|_{\mathrm{op},\infty}}
      {3\alpha(n)n^2}
 \sum_{\substack{k,l\\l>k}}|r_k\nH_{k,l}(r)-r_l\nH_{l,k}(r)|
 \leq
 \frac{\sqrt2\,C_{\mathrm{bal}}}{3n^2}
 \|D^3h\|_{\mathrm{op},\infty}.
$$
Finally, for the remainder, \eqref{eq:trace_total_rate} gives
\begin{align*}
 \left|
 \frac n{\alpha(n)}
 \sum_{\substack{k,l\\k\neq l}}r_k\nH_{k,l}(r)R_{k,l}^{(4)}(r)
 \right|
 &\leq
 \frac{\|D^4h\|_{\mathrm{op},\infty}}
      {6\alpha(n)n^3}
 \sum_{\substack{k,l\\k\neq l}}r_k\nH_{k,l}(r) \leq
 \frac{\|\bar G\|_{\mathrm{tr},\infty}}
      {12a_n^2n^2}
 \|D^4h\|_{\mathrm{op},\infty}.
\end{align*}
Here we used $a_n^2=n\alpha(n)/\beta(n)$. Consequently,
\begin{align}\label{eq:critical_pairwise_gen_comp_PRM}
\begin{aligned}
 &\left|
 \bngen(h\circ\nden)(x)
 -(\mathcal B_{a_n}^{(n)}h)(r)
 \right|\\
 &\qquad\leq
 \frac1{n^2}\left(
 \frac{\sqrt2\,C_{\mathrm{bal}}}{3}
 \|D^3h\|_{\mathrm{op},\infty}
 +\frac{\|\bar G\|_{\mathrm{tr},\infty}}{12a_*^2}
 \|D^4h\|_{\mathrm{op},\infty}
 \right).
\end{aligned}
\end{align}

Fix $t\leq T$ and $\psi\in C_b^\infty(\R^K;\R)$, and let
$\bar g_{s,t}$ be the backward function used in the proof
of Theorem \ref{thm:critical_dK_PRM}. Applying
\eqref{eq:critical_pairwise_gen_comp_PRM} in
\eqref{eq:critical_Stein_identity_PRM}, and then using
Lemma \ref{lem:critical_backward_reg_PRM} through order
four, gives
\begin{align*}
 &\left|
 \EE[\psi(\bnmu_t)]
 -\EE[\psi(\crlim^{(n,a_n)}(t))]
 \right| 
\leq
 \bar C'_{\mathrm{test}}(T)n^{-2}
 \sum_{j=1}^4\|D^j\psi\|_{\mathrm{op},\infty},
\end{align*}
where
$$
 \bar C'_{\mathrm{test}}(T)
 \dfeq T\left(
 \frac{\sqrt2\,C_{\mathrm{bal}}}{3}
 C_{\mathrm{back},3}(a_*,T)
 +\frac{\|\bar G\|_{\mathrm{tr},\infty}}{12a_*^2}
 C_{\mathrm{back},4}(a_*,T)
 \right).
$$
Since both laws have finite first moments, Lemma
\ref{lem:regularization_dK} with $\ell=4$ gives
\begin{align*}
 &\sup_{t\leq T}\wass_1\bigl(
 \operatorname{Law}(\bnmu_t),
 \operatorname{Law}(\crlim^{(n,a_n)}(t))
 \bigr) 
 \leq C'_{\mathrm{SDE},0}(T)n^{-1/2},
\end{align*}
where
$
 C'_{\mathrm{SDE},0}(T)
 \dfeq C_{\mathrm{reg}}(K,4)
 \max_{1\leq j\leq4}
 \left(\bar C'_{\mathrm{test}}(T)\right)^{1/j}.
$

The triangle inequality and the coupling estimate in the
proof of Theorem \ref{thm:critical_dK_PRM} now give
\begin{align*}
 &\sup_{t\leq T}\wass_1\bigl(
 \operatorname{Law}(\bnmu_t),
 \operatorname{Law}(\crlim^{(n,a_\infty)}(t))
 \bigr)
\leq
 C'_{\mathrm{SDE},1}(T)
 \bigl[n^{-1/2}+|a_n^{-1}-a_\infty^{-1}|\bigr],
\end{align*}
with
$
 C'_{\mathrm{SDE},1}(T)
 \dfeq C'_{\mathrm{SDE},0}(T)\vee C_{\mathrm{coup}}(T).
$

\end{proof}

\subsection{Wright–-Fisher diffusion approximation}
\label{sec:WF-diff-proof}
\begin{proof}[Proof of Theorem \ref{thm:critical_WF_dK_PRM}]
The proof follows the generator comparison and regularization
arguments in the proofs of Theorems
\ref{thm:critical_dK_PRM} and \ref{thm:critical_pairwise}.
Fix $t\leq T$, $c\in\{a_n,a_\infty\}$, and
$\psi\in C_b^\infty(\R^K;\R)$.
Let $\bar g_{s,t}$ be the backward function in Lemma
\ref{lem:critical_WF_backward_regularity_PRM} with terminal
value $\psi\circ\mathfrak r$, and put
$$
 \hat{g}_{s,t}(r)\dfeq\bar g_{s,t}(\mathfrak p(r)),
 \qquad r\in\mathcal P(\mathbb S).
$$
Thus $\hat{g}_{t,t}(r)=\psi(r)$, and the same argument as in \eqref{eq:critical_Stein_identity_PRM} gives
\begin{align*}
 &\left|
 \EE[\psi(\bnmu_t)]
 -\EE[\psi(\crlim_{\mathrm{WF}}^{(n,c)}(t))]
 \right|\\
 &\leq
 \int_0^t\EE\left[
 \left|
 \bngen(\hat{g}_{s,t}\circ\nden)(\bnX(s))
 -(\mathcal B_{\mathrm{WF},c}^{(n)}\bar g_{s,t})
   (\mathfrak p(\bnmu_s))
 \right|
 \right]\,ds.
\end{align*}
Similarly, as in the proof of the second inequality in
\eqref{eq:critical_gen_comp_PRM}, we control the contribution
of the covariance discrepancy using
$|\trace(AB)|\leq\|A\|_{\mathrm{tr}}\|B\|_{\mathrm{op}}$.
In the present setting, \eqref{eq:WF_covariance_error_PRM}
gives
\begin{align*}
 &\left\|
 a_n^{-2}\bnG(r)
 -2c^{-2}\bigl(\operatorname{diag}(r)-rr^{\mathsf T}\bigr)
 \right\|_{\mathrm{tr}} \leq
 a_*^{-2}\varepsilon_n+2|a_n^{-2}-c^{-2}|,
\end{align*}
and the same bound holds for its upper left $d\times d$
block. Consequently, the Taylor estimates used in the
preceding proofs, together with
\eqref{eq:WF_backward_regularity_PRM} and
$$
 \|D^j(\psi\circ\mathfrak r)\|_{\mathrm{op},\infty}
 \leq K^{j/2}\|D^j\psi\|_{\mathrm{op},\infty},
$$
give
\begin{align*}
 &\left|
 \EE[\psi(\bnmu_t)]
 -\EE[\psi(\crlim_{\mathrm{WF}}^{(n,c)}(t))]
 \right|\\
 &\qquad\leq C(T)\left[
 n^{-(\ell-2)}
 \sum_{j=1}^{\ell}\|D^j\psi\|_{\mathrm{op},\infty}
 +\bigl(\varepsilon_n+|a_n^{-2}-c^{-2}|\bigr)
 \sum_{j=1}^2\|D^j\psi\|_{\mathrm{op},\infty}
 \right].
\end{align*}
We take $\delta_{n,c}=\varepsilon_n+|a_n^{-2}-c^{-2}|$.
Applying Lemma~\ref{lem:regularization_dK} to the preceding
smooth-test estimate gives
\[
 \begin{aligned}
 \sup_{t\leq T}
 \wass_1\bigl(
   \operatorname{Law}(\bnmu_t),
   \operatorname{Law}(\crlim_{\mathrm{WF}}^{(n,c)}(t))
 \bigr)\leq C(T)\left[
   n^{-(\ell-2)/\ell}+\delta_{n,c}+\delta_{n,c}^{1/2}
 \right].
 \end{aligned}
\]
Since $\bnG(r)$ is nonnegative definite,
\[
 \begin{aligned}
 \varepsilon_n
 \leq\|\bar G\|_{\mathrm{tr},\infty}
   +2\sup_{r\in\mathcal P(\mathbb S)}
     \trace\bigl(\operatorname{diag}(r)-rr^\top\bigr)\leq\|\bar G\|_{\mathrm{tr},\infty}+2.
 \end{aligned}
\]
Also, $a_n,a_\infty\geq a_*$, so
$|a_n^{-2}-c^{-2}|\leq a_*^{-2}$ for
$c\in\{a_n,a_\infty\}$. Hence $0\leq\delta_{n,c}
 \leq\|\bar G\|_{\mathrm{tr},\infty}+2+a_*^{-2}$,
and therefore
 $\delta_{n,c}
 \leq
 \bigl(\|\bar G\|_{\mathrm{tr},\infty}+2+a_*^{-2}\bigr)^{1/2}
 \delta_{n,c}^{1/2}$.
Thus we can absorb this factor into $C(T)$. Using $\delta_{n,c}^{1/2}\leq\varepsilon_n^{1/2}
+|a_n^{-2}-c^{-2}|^{1/2}$ and taking $\ell=3$, we prove
the bounds in \eqref{eq:critical_WF_dK_bound_PRM} by setting
$c=a_n$ and $c=a_\infty$, respectively.
Under \eqref{eq:pairwise_balance}, taking $\ell=4$ in the
same argument improves $n^{-1/3}$ to $n^{-1/2}$ in both bounds.
The constant $C(T)$ is uniform in $n$, $t\leq T$
and $c\in\{a_n,a_\infty\}$.
\end{proof}

\section{Examples and applications}\label{sec:examples_PRM}
In this section, we study quantitative approximations for the noisy voter model with strong, moderate, and weak noise,
Moran models with mutation and selection, the heat-bath Glauber dynamics for the $p$-spin Curie--Weiss
model, and autocatalytic
reaction systems with mutations. We apply Theorems~\ref{th:ODE-approx}, \ref{thm:kantorovich_PRM},
and~\ref{thm:concentration_PRM} in the ODE regime and apply Theorem~\ref{thm:critical_dK_PRM},
\ref{thm:critical_pairwise}, and \ref{thm:critical_WF_dK_PRM} in the critical diffusion
regime. For the noisy voter model and Moran
models with selections, the scaling of the transition rates determines
which regime the models fall into. Heat-bath Glauber dynamics for the $p$-spin Curie--Weiss
model falls in the ODE regime. For
autocatalytic reactions with mutations, the balanced conversion case is in the critical regime, and the unbalanced conversion case is in the ODE regime.
We use $x\in\mathbb S^n$ to denote a configuration in this section. We retain this notation for spin systems, where $\sigma$ is usually used to denote configurations in statistical mechanics models.

\subsection{Noisy voter and neutral Moran models}
\label{subsec:voter_models_PRM}
Holley and Liggett~\cite{HolleyLiggett1975} established the foundation for the voter model, and Granovsky and Madras
\cite{GranovskyMadras1995} studied their voter model with noise.
For symmetric noise of constant order (strong noise in our case below), Cox, Peres and Steif
\cite[Theorem~1]{CoxPeresSteif2016} proved a cutoff phenomenon, while Pymar and Rivera
\cite[Theorem~1 and Proposition~8]{PymarRivera2024}
studied the density fluctuations under stationarity.
In the weak-noise regime, Aljovin, Jara and Xiang~\cite{AljovinJaraXiang2024}
studied thermalisation and quantitative approximation to a nonlinear diffusion. Their Theorem~3.1
gives a Wright--Fisher approximation for the density at
time $nt$ in Kantorovich distance, with error at most
$C(1\wedge t^{1/4})n^{-1/2}$.
The fixed-interval rate follows from
Theorem~\ref{thm:critical_WF_dK_PRM}.
Keeping the derivative decay in its generator comparison
proof gives the stated dependence on time.

Under moderate or strong noise regimes, the noisy voter model falls into an ODE regime.
Theorems~\ref{th:ODE-approx} and~\ref{thm:concentration_PRM} give
moment and concentration bounds,
while Theorem~\ref{thm:kantorovich_PRM} gives a quantitative bound for the central limit approximation of the density fluctuations.
To our knowledge, no previous work gives a quantitative bound of this form covering all noise regimes. 

\subsubsection{Dynamics and empirical characteristics}
\label{subsubsec:nvm_common_computation_PRM}

We relabel the states as $\mathbb S=\{0,1\}$, so $\mu_t^{(n)}(1)$ is
the density of type-$1$ sites. Fix $a,b>0$, and a sequence $\gamma_n>0$.
We consider a deterministic initial empirical measure $\mu_0^{(n)}$.
Each site chooses a site uniformly at random and copies its state, and it may choose itself. The dynamics has unit
copying rate and independent spontaneous flips at rates $a\gamma_n/n$
from $0$ to $1$ and $b\gamma_n/n$ from $1$ to $0$. The per-site flip
rates at empirical state $r=(r_0,r_1)$ are 
\begin{equation}
\label{eq:nvm_common_H_PRM}
 H_{0,1}^{(n)}(r)=r_1+\frac{a\gamma_n}{n},
 \qquad H_{1,0}^{(n)}(r)=r_0+\frac{b\gamma_n}{n}.
\end{equation}
Using $\bar{\canvec}_{1,0}=-\bar{\canvec}_{0,1}$
in \eqref{eq:cov-def-comp}, we obtain
\begin{align*}
 F^{(n)}(r)
 =\frac{\gamma_n}{n}[a-(a+b)r_1]\bar{\canvec}_{0,1}, \ 
 G^{(n)}(r)
 =\left[2r_1(1-r_1)+\frac{\gamma_n}{n}\{a(1-r_1)+br_1\}\right]
   \bar{\canvec}_{0,1}\bar{\canvec}_{0,1}^{\top}.
\end{align*}

\subsubsection{Identification with a neutral Moran model}
\label{subsubsec:moran_identification_PRM}
We illustrate that a neutral Moran model with independent mutation can be represented
as a noisy voter model. The classical model can be found in Moran~\cite{Moran1958}.
For reproduction, each individual reproduces at rate one, and its offspring inherits its type and replaces a uniformly chosen individual. For mutations, each individual mutates independently, with mutations $0\to1$ and $1\to0$ occurring at rates $a\gamma_n/n$ and $b\gamma_n/n$, respectively.

We write $r_1$ for the type-$1$ density, so the type-$0$ density is $r_0=1-r_1$. For a fixed type-$0$ individual. There are $nr_1$ type-$1$ parents, and each parent's offspring replaces this particular individual
with probability $1/n$. As a result, the contribution of reproduction to its rate $0\to1$ is therefore $(nr_1)/n=r_1$. Similarly, the reproduction contribution to a type-$1$ individual's $1\to0$ rate is $1-r_1$. Including mutation, the per-individual flip rates are
\begin{equation*}
 H_{0,1}^{(n)}(r)=r_1+\frac{a\gamma_n}{n},
 \qquad
 H_{1,0}^{(n)}(r)=1-r_1+\frac{b\gamma_n}{n}.
\end{equation*}
These are the rates of the noisy-voter model in
\eqref{eq:nvm_common_H_PRM}.

\subsubsection{Moderate or strong noise in the ODE regime}
\label{subsec:nvm_moderate_noise_PRM}

Assume that $\gamma_n\to\infty$ and $\sup_n\gamma_n/n<\infty$.
We call the noise \emph{moderate} when $\gamma_n/n\to0$ and \emph{strong} when
$\gamma_n/n\to\theta\in(0,\infty)$. For example, $a=b=1/2$ and
$\gamma_n=n$ give the flip rate of
\cite[Section~1]{CoxPeresSteif2016}. Take $\alpha(n)=\gamma_n/n$
and $\beta(n)=1$. Then $\bnmu_t=\mu_{nt/\gamma_n}^{(n)}$,
$a_n=\sqrt{\gamma_n}$, and
\[
 \bar F^{(n)}(r)=[a-(a+b)r_1]\bar{\canvec}_{0,1},\qquad \bar G^{(n)}(r)=G^{(n)}(r).
\]
Thus Assumptions~\ref{assum:scaling-1} and~\ref{assum:ode-regime} hold. 
The solution of \eqref{eq:ODE-approx} with $\rho_0^{(n)}=\mu_0^{(n)}$
has type-$1$ coordinate
\[
 \rho_t^{(n)}(1)=\frac{a}{a+b}
 +\left(\mu_0^{(n)}(1)-\frac{a}{a+b}\right)e^{-(a+b)t}.
\]

Theorems~\ref{th:ODE-approx} and~\ref{thm:kantorovich_PRM} apply with
$a_n=\sqrt{\gamma_n}$. The first gives  the $p$th moment bound of $\sup_{t\leq T}\|\bnmu_t-\rho_t^{(n)}\|$. The second bounds the $\wass_1$ distance
between the law of $\fluct_t^{(n)}$ and the time-$t$ law of the
linear diffusion in \eqref{eq:Y_PRM_SDE}.
Theorem~\ref{thm:concentration_PRM} applies with $\|D\bar F\|_{\mathrm{op},\infty}=a+b$,
$\delta_T=\delta e^{-(a+b)T}$, and
 $\|\bar G\|_{\mathrm{tr},\infty}
 \leq 1+2\max\{a,b\}\sup_{n\geq1}\frac{\gamma_n}{n}$.
The constant and exponent scale in \eqref{eq:concentration_main_PRM} are
\[
 C_{\mathrm{conc}}(T,\delta)
 =\frac{\delta^2e^{-2(a+b)T}}
 {2\bigl(T\|\bar G\|_{\mathrm{tr},\infty}
        +\sqrt2\,\delta e^{-(a+b)T}/3\bigr)},
 \qquad a_n^2\wedge n=\gamma_n\wedge n\asymp\gamma_n.
\]

\subsubsection{Weak noise critical regime: Wright--Fisher approximation}
\label{subsec:nvm_weak_noise_PRM}
Assume $\gamma_n\to\gamma_\infty\in(0,\infty)$. We call this noise intensity the
\emph{weak-noise regime}. Under this regime, we have $\gamma_n/n\to0$, and spontaneous flip rates are both of order $n^{-1}$. Now take
$\alpha(n)=n^{-1}$ and $\beta(n)=1$, so $a_n^{-2}=1$. We speed up time by a factor $n$, and we observe the process by $\bnmu_t=\mu_{nt}^{(n)}$. Assumption~\ref{assum:critical_regime_PRM} \ref{item:FF-exten} fails at the boundary of densities, so neither
Theorem~\ref{thm:critical_dK_PRM} nor Theorem \ref{thm:critical_pairwise}
applies. However, Lemma~\ref{lem:critical_WF_backward_regularity_PRM}
provides the backward regularity estimate for the Wright--Fisher diffusion.
Using this estimate in the same generator comparison argument, Theorem~\ref{thm:critical_WF_dK_PRM} gives a quantitative Wright--Fisher approximation in the critical regime.

Let $\crlim_{\mathrm{WF}}^{(n,1)}
=(\crlim_{\mathrm{WF}}^{(n,1)}(\cdot,0),
  \crlim_{\mathrm{WF}}^{(n,1)}(\cdot,1))$
be the two-type Wright--Fisher diffusion with
$\crlim_{\mathrm{WF}}^{(n,1)}(0)=\mu_0^{(n)}$.
We consider the type-$1$ density satisfies
\[
 \begin{aligned}
 d\crlim_{\mathrm{WF}}^{(n,1)}(t,1)
 =\gamma_n\bigl[a-(a+b)\crlim_{\mathrm{WF}}^{(n,1)}(t,1)\bigr]\,dt+\sqrt{2\crlim_{\mathrm{WF}}^{(n,1)}(t,1)
                   \bigl(1-\crlim_{\mathrm{WF}}^{(n,1)}(t,1)\bigr)}\,dW(t),
 \end{aligned}
\]
where $W$ is a standard Brownian motion. Note that
$\crlim_{\mathrm{WF}}^{(n,1)}(t,0)=1-\crlim_{\mathrm{WF}}^{(n,1)}(t,1)$,
and the covariance, 
\[
 \begin{aligned}
 \bar G^{(n)}(r)
 =2\bigl(\operatorname{diag}(r)-rr^\top\bigr)+\frac{\gamma_n}{n}[a(1-r_1)+br_1]\,
       \bar{\canvec}_{0,1}\bar{\canvec}_{0,1}^{\top}.
 \end{aligned}
\]
Thus diffusion contains the exact drift and the leading term in the covariance. 
Since $\|\bar{\canvec}_{0,1}\|^2=2$, the covariance discrepancy
in \eqref{eq:WF_covariance_error_PRM} is
\[
 \varepsilon_n
 =\frac{2\gamma_n}{n}\max\{a,b\}
 \leq\frac{2\max\{a,b\}}{n}\sup_{m\geq1}\gamma_m.
\]
Therefore, for every $T>0$, Theorem~\ref{thm:critical_WF_dK_PRM} gives, 
\[
 \sup_{0\leq t\leq T}
 \wass_1\bigl(\operatorname{Law}(\bnmu_t),
             \operatorname{Law}(\crlim_{\mathrm{WF}}^{(n,1)}(t))\bigr)
 \leq C_Tn^{-1/2}.
\]

However, we can keep track of the time variable, instead of bounding it directly by a constant depending on $T$. By keeping the derivative decay in
\eqref{eq:WF_neutral_derivative_decay_PRM}, the argument in \cite[proof of Theorem~3.1]{AljovinJaraXiang2024} further gives a time dependence bound.
For every $t\geq0$,
\[
\wass_1\bigl(\operatorname{Law}(\bnmu_t(1)),
             \operatorname{Law}(\crlim_{\mathrm{WF}}^{(n,1)}(t,1))\bigr)
 \leq C(1\wedge t^{1/4})n^{-1/2}, 
\]
where $C$ depends only on $a,b$ and $\sup_n\gamma_n$. When $\gamma_n\equiv\gamma>0$, this estimate follows from
\cite[Theorem~3.1]{AljovinJaraXiang2024} with $a\gamma$ and $b\gamma$
in place of their mutation parameters $a$ and $b$, respectively.

\subsection{Moran reproduction with mutations and selections}
\label{subsec:selected_moran_PRM}
Unlike the neutral Moran model in Section \ref{subsubsec:moran_identification_PRM}, when reproduction rates differ between types, the model is not the noisy-voter type model. It admits an interpretation as biased copying with mutations. 
For two-type Moran models with mutations and selections,
Cordero~\cite[Proposition~3.1 and Theorem~3.4]{Cordero2017}
proves a deterministic limit and a central limit theorem. For multiple types, Cloez and
Corujo~\cite[Theorems~1.3, 1.5 and 1.8]{CloezCorujo2022}
establish quantitative propagation of chaos and Gaussian
fluctuations, with uniform-in-time moment bounds. 

We consider two regimes for the Moran model. Each individual reproduces
at rate $1$ plus a type-dependent rate proportional to $\gamma_n/n$.
The mutation rates are proportional to $\gamma_n/n$.
If $\gamma_n\to\infty$ and $\sup_n(\gamma_n/n)<\infty$, the covariance
rate on the time scale $t\mapsto nt/\gamma_n$ is $O(\gamma_n^{-1})$.
The model therefore falls in our ODE regime. If
$\gamma_n\to\gamma_\infty\in(0,\infty)$, the covariance on the time
scale $t\mapsto nt$ has a nonzero limit.
In this critical regime, we approximate the empirical type
frequencies by a multitype Wright--Fisher diffusion with selections
and mutations.

Fix $\mathbb S=\{1,\ldots,K\}$, selection parameters $s_k\geq0$,
mutation parameters $\vartheta_k>0$, a sequence $\gamma_n>0$, and an initial empirical measure $\mu_0^{(n)}$. A type-$k$
individual reproduces at rate $1+s_k\gamma_n/n$. For reproduction with selections, its offspring inherits
its type and replaces a uniformly chosen individual, possibly its parent.
For independent mutations, each type-$k$ individual mutates to type $l\ne k$ at rate
$\vartheta_l\gamma_n/n$. Its transition rate to type $l$ is therefore
\begin{equation}
\label{eq:selected_moran_rates_PRM}
 H_{k,l}^{(n)}(r)=\left(1+\frac{s_l\gamma_n}{n}\right)r_l
                    +\frac{\vartheta_l\gamma_n}{n},\qquad k\ne l.
\end{equation}

We write $\bar s(r)=\sum_l s_lr_l$ and $\Theta=\sum_l\vartheta_l$.
Thus
\eqref{eq:cov-def-comp} gives
\begin{align*}
 &F^{(n)}(k,r)=\frac{\gamma_n}{n}
 \bigl[r_k(s_k-\bar s(r))+\vartheta_k-\Theta r_k\bigr],\\
 &G^{(n)}(r)=\sum_{k<l}\left[2r_kr_l+\frac{\gamma_n}{n}
 \bigl\{(s_k+s_l)r_kr_l+\vartheta_lr_k+\vartheta_kr_l\bigr\}\right]
 \bar{\canvec}_{k,l}\bar{\canvec}_{k,l}^{\top}.
\end{align*}

\subsubsection{Moderate or strong noise in the ODE regime}
\label{subsubsec:selected_moran_ode_regime_PRM}

Assume $\gamma_n\to\infty$ and $\sup_n(\gamma_n/n)<\infty$.
As in Section~\ref{subsec:nvm_moderate_noise_PRM}, we call the
scaling moderate when $\gamma_n/n\to0$ and strong when
$\gamma_n/n\to\theta\in(0,\infty)$. In both cases, take
$\alpha(n)=\gamma_n/n$ and $\beta(n)=1$. Then
\[
 \bnmu_t=\mu_{nt/\gamma_n}^{(n)},\qquad
 a_n=\sqrt{\gamma_n},\qquad
 \bar F^{(n)}=\frac{n}{\gamma_n}F^{(n)},\qquad
 \bar G^{(n)}=G^{(n)}.
\]
Since $a_n\to \infty$, the model falls into the ODE
regime. Assumptions~\ref{assum:scaling-1} and~\ref{assum:ode-regime} hold. 
We take the initial condition as $\rho_0^{(n)}=\mu_0^{(n)}$, and the ODE in \eqref{eq:ODE-approx}
becomes
\[
 \dot\rho_t^{(n)}(k)
 =\rho_t^{(n)}(k)\bigl[s_k-\bar s(\rho_t^{(n)})\bigr]
  +\vartheta_k-\Theta\rho_t^{(n)}(k),
 \qquad k\in\mathbb S.
\]

On a fixed interval $[0,T]$, Theorem~\ref{th:ODE-approx} can be applied to give the $p$th moment of $\sup_{t\leq T}\|\bnmu_t-\rho_t^{(n)}\|$. For the fluctuation
$\fluct_t^{(n)}=\sqrt{\gamma_n}(\bnmu_t-\rho_t^{(n)})$,
Theorem~\ref{thm:kantorovich_PRM} gives
\[
 \sup_{0\leq t\leq T}
 \wass_1\bigl(\operatorname{Law}(\fluct_t^{(n)}),
              \operatorname{Law}(Y_{0,0}^{(n)}(t))\bigr)
 \leq C_T\max\{\gamma_n^{-1/4},\,\gamma_n^{1/6}n^{-1/3}\},
\]
where $Y_{0,0}^{(n)}$ is the linear diffusion in
\eqref{eq:Y_PRM_SDE}.

We can apply Theorem~\ref{thm:concentration_PRM} with 
$\|D\bar F\|_{\mathrm{op},\infty}\leq \Theta+\max_k s_k+(\sum_k s_k^2)^{1/2}$ and
$\delta_T=\delta e^{-\|D\bar F\|_{\mathrm{op},\infty}T}$. The covariance trace satisfies
\[
 \|\bar G\|_{\mathrm{tr},\infty}
 =2\sup_{n,r}\sum_{k\ne l}r_kH_{k,l}^{(n)}(r)
 \leq 2\left[1+(\max_k s_k+\Theta)\sup_n\frac{\gamma_n}{n}\right].
\]
Thus \eqref{eq:concentration_main_PRM} gives the upper bound
$2\exp\{-C_{\mathrm{conc}}(T,\delta)(\gamma_n\wedge n)\}$, where
\[
 C_{\mathrm{conc}}(T,\delta)
 =\frac{\delta_T^2}
        {2T\|\bar G\|_{\mathrm{tr},\infty}+(2\sqrt2/3)\delta_T},
 \qquad \gamma_n\wedge n\asymp\gamma_n.
\]

In the strong-noise case, take $\gamma_n=n$.
For $K=2$ and $s_1=0$, this fixed-parameter model agrees with Cordero's model after relabeling  in \cite[Section~2]{Cordero2017}. Identify our type $2$ with his
type $0$ and take $N=n$, $s_N=s_2$, $u_N=\Theta$,
$\nu_0=\vartheta_2/\Theta$, and $\nu_1=\vartheta_1/\Theta$.
It follows that the rates are 
\[
 H_{1,2}^{(n)}(r)=(1+s_2)r_2+\vartheta_2,\qquad
 H_{2,1}^{(n)}(r)=r_1+\vartheta_1.
\]
If we multiply the rates above by $nr_1$ and $nr_2$, respectively, then we obtain the birth
and death rates following equation~(2.1) in
\cite[Section~2]{Cordero2017}.

For general $K$, we take $N=n$, $E=\mathbb S$, and $\eta(k)=nr_k$
in \cite[Section~1]{CloezCorujo2022}. We define $Q_{k,l}=\vartheta_l$
for $k\ne l$, $Q_{k,k}=-(\Theta-\vartheta_k)$, and
$V_\mu(k,l)=1+s_l$ for every $\mu$. The total rate of replacing
a type-$k$ individual by type $l$ is then
\[
 \eta(k)\left[Q_{k,l}+\frac{\eta(l)}N V_\mu(k,l)\right]
 =nr_k\bigl[\vartheta_l+(1+s_l)r_l\bigr]
 =nr_kH_{k,l}^{(n)}(r),\qquad k\ne l.
\]
Thus the two count processes have the same generator. Therefore, our results apply to their models.

\subsubsection{The weak noise case:
\texorpdfstring{$\gamma_n\to\gamma_\infty\in(0,\infty)$}
{gamma(n) tending to a finite positive limit}}

Under the weak noise regime, we take $\alpha(n)=n^{-1}$
and $\beta(n)=1$, so $a_n=a_\infty=1$.
We observe the process on the time scale $nt$ and write
$\bnmu_t=\mu_{nt}^{(n)}$. On this scale, the normalized
drift and covariance are $\bnF=nF^{(n)}$ and $\bnG=G^{(n)}$.
Assumption~\ref{assum:critical_regime_PRM} \ref{item:FF-exten} fails at the
boundary. We instead apply
Theorem~\ref{thm:critical_WF_dK_PRM}.

The drift coordinates have expression
\[
 \bnF(k,r)=\gamma_n
 \bigl[r_k(s_k-\bar s(r))+\vartheta_k-\Theta r_k\bigr].
\]
We also have $\sup_n\gamma_n<\infty$.
Thus \eqref{eq:WF_drift_regularity_PRM} holds.

Let $\crlim_{\mathrm{WF}}^{(n,1)}$ be the
diffusion in \eqref{eq:WF_diffusion_SDE_PRM},
with $c=1$. We can write it as
\[
 \begin{aligned}
 d\crlim_{\mathrm{WF}}^{(n,1)}(t)
 =\bnF\bigl(\crlim_{\mathrm{WF}}^{(n,1)}(t)\bigr)\,dt+\sum_{k<l}
 \sqrt{2\crlim_{\mathrm{WF}}^{(n,1)}(t,k)
          \crlim_{\mathrm{WF}}^{(n,1)}(t,l)}\,
 \bar{\canvec}_{k,l}\,dW^{k,l}(t),
 \end{aligned}
\]
with $\crlim_{\mathrm{WF}}^{(n,1)}(0)=\mu_0^{(n)}$
and independent standard Brownian motions $W^{k,l}$.

The covariance discrepancy in
\eqref{eq:WF_covariance_error_PRM} satisfies
\[
 \varepsilon_n
 =\frac{2\gamma_n}{n}
   \sup_{r\in\mathcal P(\mathbb S)}
   \sum_{k=1}^K(s_kr_k+\vartheta_k)(1-r_k)
 \leq\frac{2\gamma_n}{n}(\max_k s_k+\Theta).
\]
Moreover, \eqref{eq:pairwise_balance} holds with
$C_{\mathrm{bal}}=(\max_k s_k+\Theta)\sup_n\gamma_n$. Therefore, Theorem~\ref{thm:critical_WF_dK_PRM} gives,
for every fixed $T>0$,
\begin{equation}
\label{eq:selected_moran_weak_dK_PRM}
 \begin{aligned}
 \sup_{0\leq t\leq T}\wass_1\bigl(
   \operatorname{Law}(\bnmu_t),
   \operatorname{Law}(\crlim_{\mathrm{WF}}^{(n,1)}(t))
 \bigr)\leq C'_{\mathrm{WF},0}(T)
       \bigl[n^{-1/2}+\varepsilon_n^{1/2}\bigr]
 \leq C_Tn^{-1/2},
 \end{aligned}
\end{equation}
where $C_T$ depends only on $K,T,\max_k s_k,\Theta$ and
$\sup_n\gamma_n$, and is independent of $n$ and the initial
empirical state.

\subsection{Continuous-time Glauber dynamics for the \texorpdfstring{$p$}{p}-spin Curie--Weiss model}
\label{subsec:p_spin_glauber_PRM}
The Curie--Weiss model is a mean-field spin system, which is defined as the Ising model on complete graphs.
Ding, Lubetzky and Peres~\cite{DingLubetzkyPeres2009} and Levin, Luczak and Peres~\cite{LevinLuczakPeres2010} investigated the mixing properties in different temperature regimes, including studies on cutoff phenomenon, critical power law,
and metastability. Mukherjee, Son and Bhattacharya~\cite{MukherjeeSonBhattacharya2021}
studied magnetization fluctuations in the $p$-spin model,
including Gaussian and non-Gaussian limits. Samanta, Mukherjee and
Zhang~\cite{SamantaMukherjeeZhang2026} investigated mixing properties
and metastability in the $p$-spin Curie--Weiss model.
Here we study the quantitative approximations for heat-bath Glauber dynamics for the $p$-spin Curie--Weiss model.

\subsubsection{Dynamics and empirical characteristics}
\label{subsubsec:p_spin_model_PRM}
Let $\mathbb S=\{-1,+1\}$ be the set of spins and fix an integer $p\geq2$. We take the inverse
temperature $\beta\geq0$, and external field $h\in\mathbb R$.
The magnetization of $x\in\mathbb S^n$ is defined as
$m(x)=n^{-1}\sum_{i=1}^nx_i$. Due to the mean-field setting, the Gibbs measure and its
partition function are
\[
 \pi_n(x)=Z_n^{-1}e^{n[\beta m(x)^p+hm(x)]},
 \qquad Z_n=\sum_{x\in\mathbb S^n}e^{n[\beta m(x)^p+hm(x)]}.
\]
Each site is updated at rate one by sampling from its conditional
Gibbs measure given the other spins. We write $x^{(i)}$ for
the configuration obtained by flipping spin $i$. The transition
$x\to x^{(i)}$ has rate
$\pi_n(x^{(i)})/[\pi_n(x)+\pi_n(x^{(i)})]$. For the Glauber dynamics, $\pi_n$
is the unique invariant law.

At empirical state $r$, we take the magnetization $m=m(r)=r(+1)-r(-1)$.
If site $i$ currently has spin $x_i=s$, the other spins contribute to the magnetization by
$z=n^{-1}\sum_{j\ne i}x_j=m(x)-s/n$.
By straightforward computations, the conditional probability of assigning $+1$ is
$(1+e^{-2\theta_n(z)})^{-1}=(1+\tanh\theta_n(z))/2$, where
\[
 \theta_n(z)=h+\frac{n\beta}{2}
 \bigl[(z+1/n)^p-(z-1/n)^p\bigr].
\]
Since each site is updated at rate one, its flip rate is the
conditional probability of assigning $-s$ when its current spin is $s$.
We write $\lambda_{n,s}(m)=\tanh\theta_n(m-s/n)$, and it follows that
\[
 H_{s,-s}^{(n)}(r)=\frac{1-s\lambda_{n,s}(m)}{2},
 \qquad s\in\{-1,+1\}.
\]

We have that $r(\pm1)=(1\pm m)/2$ and a flip changes the magnetization
by $\pm2/n$. Therefore, we obtain
\begin{align*}
 F_n(m)&=\frac{1-m}{2}\bigl(1+\lambda_{n,-1}(m)\bigr)
         -\frac{1+m}{2}\bigl(1-\lambda_{n,+1}(m)\bigr),\\
 G_n(m)&=(1-m)\bigl(1+\lambda_{n,-1}(m)\bigr)
         +(1+m)\bigl(1-\lambda_{n,+1}(m)\bigr).
\end{align*}

\subsubsection{ODE regime and quantitative approximations}
\label{subsubsec:p_spin_ode_PRM}
Under the ODE regime, we take
$\alpha(n)=\beta(n)=1$, and we have 
$a_n=\sqrt n$. We do not speed up time, and we take $\bnmu_t=\mu_t^{(n)}$. Here $\beta(n)$ is the covariance normalization and
$\beta$ is the inverse temperature. Thus
\eqref{eq:cov-def-comp} gives
\[
 \bar F^{(n)}(r)=\frac{F_n(m)}2\,\bar{\canvec}_{-1,+1},
 \qquad
 \bar G^{(n)}(r)=\frac{G_n(m)}4\,
             \bar{\canvec}_{-1,+1}\bar{\canvec}_{-1,+1}^{\top}.
\]

We fix an initial empirical measure $\mu_0^{(n)}$. Let
$\rho^{(n)}$ solve \eqref{eq:ODE-approx} with
$\rho_0^{(n)}=\mu_0^{(n)}$. Its magnetization satisfies the mean-field equation
\[
 \frac{d}{dt}m(\rho_t^{(n)})=F_n\bigl(m(\rho_t^{(n)})\bigr),
 \qquad m(\rho_0^{(n)})=m(\mu_0^{(n)}).
\]
The formula above for $\theta_n$ gives
$F_n(m)\to\tanh(h+p\beta m^{p-1})-m$ as $n\to\infty$,  uniformly for $m\in[-1,1]$. We have that Assumptions~\ref{assum:scaling-1} and~\ref{assum:ode-regime} hold. 

On each fixed interval $[0,T]$, Theorem~\ref{th:ODE-approx} bounds
the $q$th moment of $\sup_{t\leq T}\|\bnmu_t-\rho_t^{(n)}\|$. Theorem~\ref{thm:kantorovich_PRM} gives a quantitative bound $C_Tn^{-1/6}$ in $\wass_1$ between the law of the fluctuation
$\fluct_t^{(n)}=\sqrt n(\bnmu_t-\rho_t^{(n)})$ and the time-$t$
law of the linear diffusion in \eqref{eq:Y_PRM_SDE}.

Theorem~\ref{thm:concentration_PRM} applies with
$\|D\bar F\|_{\mathrm{op},\infty}\leq 2+\beta p(p-1)3^{p-2}$ and $\delta_T=\delta e^{-\|D\bar F\|_{\mathrm{op},\infty}T}$.
Since $\|\bar G\|_{\mathrm{tr},\infty}
=\tfrac12\sup_{n\geq1,\,|m|\leq1}G_n(m)\leq2$,
\eqref{eq:concentration_main_PRM} gives the bound
$2e^{-C_{\mathrm{conc}}(T,\delta)n}$, where
\[
 C_{\mathrm{conc}}(T,\delta)=\frac{\delta_T^2}
 {2T\|\bar G\|_{\mathrm{tr},\infty}+(2\sqrt2/3)\delta_T}.
\]
For fixed $p\geq 2$, the inverse temperature $\beta\geq 0$, the external magnetic field $h\in \mathbb R$ and time horizon $T$, the
constants are independent of $n$ and the initial
empirical state.

\subsection{Autocatalytic reactions with mutations}
\label{subsec:balanced_conversion_PRM}

Togashi and Kaneko~\cite{TogashiKaneko2001} studied a cycle of autocatalytic reactions. In such reactions each species promotes the conversion of other species into its own type. We write $S_{K+1}=S_1$. The
network has reactions as follows: 
\[
 S_i+S_{i+1}\longrightarrow 2S_{i+1},\qquad
 \varnothing\rightleftarrows S_i,\qquad 1\leq i\leq K.
\]
In a recent work \cite{FuEtAl2026Fragility}, the authors established  a rigorous stochastic averaging principle for this model, which allows one to describe the dynamics of the slow variable (total population size) in terms of ergodic moments of the fast species. See \cite[Theorem 3.1]{FuEtAl2026Fragility}.
For $K=2$, KhudaBukhsh and Xiang
\cite{KhudaBukhshXiang2026}
derived cutoff profiles for a fixed-population model with
nonlinear kinetics and weak mutations. Boundary initial
conditions introduce an additional logarithmic growth phase, changing both the leading mixing time and the profile, from a Gaussian translation profile in the bulk initial condition to a mutation-dependent Gaussian-mixture profile.

We distinguish two kinetic laws for the same fixed-population
reaction network. In the first case, opposite conversion fluxes cancel
identically, which makes the model fall into the critical regime. In the second case, conversion produces a nonzero
leading drift, which instead makes the model fall into the ODE regime. 

We consider $K$ species and define the set of species as  $\mathbb S=\{1,\ldots,K\}$, for $K\geq2$. We write $r_i=z_i/n$ for
the density of species $i$, where $\sum_i z_i=n$. Let
$\mathcal E$ be the set of the directed edges containing each directed interaction and its
opposite. We furthermore assume that the undirected graph is connected. The conversion
reactions are $S_i+S_j\longrightarrow2S_j$ for $i\to j\in\mathcal E$. Take $c_{i\to j}>0$ on $\mathcal E$ and set
$c_{i\to j}=0$ otherwise. Mutation
$S_i\longrightarrow S_j$ is allowed for every $i\ne j$, with
mutation parameters $\vartheta_j>0$. We define $\Theta=\sum_j\vartheta_j$.
The count
chain is irreducible on $\{z\in\mathbb N_0^K:\sum_i z_i=n\}$.
We take $\mu_0^{(n)}$ as an initial condition.

\subsubsection{Balanced conversion in the critical regime}
\label{subsubsec:regulated_TK_critical_PRM}

Under the balanced conversion case, we assume that
$c_{i\to j}=c_{j\to i}$. The transition rate is given by
\[
 H_{i,j}^{(n)}(r)=c_{i\to j}r_ir_j^2+\frac{\vartheta_j}{n},
 \qquad i\ne j.
\]
From the transition rate, we have that the total conversion rate
decays quadratically when either species becomes rare. We write
$\vartheta=(\vartheta_1,\ldots,\vartheta_K)^\top$.
By taking the sum, the total mutation rate is
$\sum_j\vartheta_j(1-r_j)=\Theta-\sum_j\vartheta_jr_j$.

In this case, we take $\alpha(n)=n^{-1}$ and $\beta(n)=1$, so
$a_n=a_\infty=1$. We observe the process on the time scale
$nt$ and write $\bnmu_t=\mu_{nt}^{(n)}$.
We furthermore define $\crdiff_{\mathrm{AC}}(r)$ by its columns
\[
 \bigl(\crdiff_{\mathrm{AC}}(r)\bigr)_{\cdot,i\to j}
 =\sqrt{c_{i\to j}}\,r_ir_j\bar{\canvec}_{i,j},
 \qquad i\to j\in\mathcal E.
\]
By \eqref{eq:cov-def-comp}, we have the normalized drift and covariance
\[
 \begin{aligned}
 \bnF(r)=\vartheta-\Theta r, \qquad
 \bnG(r)=\crdiff_{\mathrm{AC}}(r)\crdiff_{\mathrm{AC}}(r)^\top
 +\frac1n\sum_{i<j}(\vartheta_jr_i+\vartheta_ir_j)
       \bar{\canvec}_{i,j}\bar{\canvec}_{i,j}^{\top}.
 \end{aligned}
\]

In the diffusion, we retain the exact drift and the full
covariance. Write
\[
 Q_\vartheta(r)=\sum_{i<j}(\vartheta_jr_i+\vartheta_ir_j)
 \bar{\canvec}_{i,j}\bar{\canvec}_{i,j}^{\top},
 \qquad \vartheta_{\min}=\min_i\vartheta_i>0.
\]
It follows that
$\bnG=\crdiff_{\mathrm{AC}}\crdiff_{\mathrm{AC}}^\top
+n^{-1}Q_\vartheta$. We define
 $S_\vartheta(r)
 =\left(Q_\vartheta(r)+\frac{\mathbf1\mathbf1^\top}{K}\right)^{1/2}
  -\frac{\mathbf1\mathbf1^\top}{K}$.
This factor is smooth in a neighborhood of the simplex. We also have that $S_\vartheta(r)S_\vartheta(r)^\top=Q_\vartheta(r)$. We take $d_{\mathrm{noise}}=|\mathcal E|+K$.
Define $\crdiff^{(n)}(r)$ as the
$K\times d_{\mathrm{noise}}$ matrix whose first
$|\mathcal E|$ columns are those of
$\crdiff_{\mathrm{AC}}(r)$ and whose remaining $K$ columns
are those of $n^{-1/2}S_\vartheta(r)$. On the simplex,
\[
 \begin{aligned}
 \crdiff^{(n)}(r)\crdiff^{(n)}(r)^\top
 =\crdiff_{\mathrm{AC}}(r)\crdiff_{\mathrm{AC}}(r)^\top
   +\frac1n S_\vartheta(r)S_\vartheta(r)^\top=\bnG(r).
 \end{aligned}
\]

\begin{remark}
Here $d_{\mathrm{noise}}$ is independent of $n$ and counts
the Brownian coordinates for the conversion and mutation
noises. Although Theorems~\ref{thm:critical_dK_PRM} and
\ref{thm:critical_pairwise} are stated for square factors,
their proofs also apply to $K\times d_{\mathrm{noise}}$
factors driven by a $d_{\mathrm{noise}}$-dimensional
Brownian motion. Indeed, the matrix estimate used in
the proofs remains valid:
\[
 \trace(MM^\top)\leq K\|M\|_{\mathrm{op}}^2,
 \qquad M\in\R^{K\times d_{\mathrm{noise}}}.
\]
Thus both theorems remain valid for these rectangular
factors under the same uniform regularity assumptions,
with the same quantitative bounds.
\end{remark}

For $z\in\mathcal P(\mathbb S)$, let
$\crlim_{\mathrm{AC},z}^{(n,1)}$ denote the diffusion in
\eqref{eq:critical_diffusion_SDE_PRM} with $c=1$, and we have that
\[
 \begin{aligned}
 d\crlim_{\mathrm{AC},z}^{(n,1)}(t)
 =\bnF\bigl(\crlim_{\mathrm{AC},z}^{(n,1)}(t)\bigr)\,dt
 +\crdiff_{\mathrm{AC}}
       \bigl(\crlim_{\mathrm{AC},z}^{(n,1)}(t)\bigr)\,dW(t)+n^{-1/2}S_\vartheta
       \bigl(\crlim_{\mathrm{AC},z}^{(n,1)}(t)\bigr)\,dB(t),
 \end{aligned}
\]
with initial condition $\crlim_{\mathrm{AC},z}^{(n,1)}(0)=z$. $W=(W^{i\to j})_{i\to j\in\mathcal E}$ is an
$|\mathcal E|$-dimensional standard Brownian motion, with
independent coordinates indexed by the directed edges in
$\mathcal E$. $B$ is an
independent standard $K$-dimensional Brownian motion.
We have $\bnG_{ii}(r)=\vartheta_i/n>0$ when $r_i=0$.
Since opposite conversion fluxes are equal, we obtain
\[
 \sum_{i<j}
 \left|r_iH_{i,j}^{(n)}(r)-r_jH_{j,i}^{(n)}(r)\right|
 =\frac1n\sum_{i<j}|\vartheta_jr_i-\vartheta_ir_j|
 \leq\frac{\Theta}{n}.
\]
Thus \eqref{eq:pairwise_balance} holds with
$C_{\mathrm{bal}}=\Theta$.
Theorem~\ref{thm:critical_pairwise} improves the bound in Theorem~\ref{thm:critical_dK_PRM} to
\[
 \sup_{0\leq t\leq T}
 \wass_1\bigl(
   \operatorname{Law}(\bnmu_t),
   \operatorname{Law}(
     \crlim_{\mathrm{AC},\mu_0^{(n)}}^{(n,1)}(t))
 \bigr)
 \leq C_{\mathrm{AC},0}(T)n^{-1/2}.
\]
The constant depends
only on $K,T,\Theta,\vartheta_{\min}^{-1}$ and
$\max_{i\to j\in\mathcal E}c_{i\to j}$.
It is independent of $n$ and the initial empirical state.

\begin{remark}[Comparison using the leading covariance]
In the comparison diffusion, we may also retain only the leading
conversion covariance. For $z\in\mathcal P(\mathbb S)$, let
$\widetilde{\crlim}_{\mathrm{AC},z}^{(n,1)}$ solve
\[
 \begin{aligned}
 d\widetilde{\crlim}_{\mathrm{AC},z}^{(n,1)}(t)
=\bnF\bigl(\widetilde{\crlim}_{\mathrm{AC},z}^{(n,1)}(t)\bigr)\,dt
+\crdiff_{\mathrm{AC}}
    \bigl(\widetilde{\crlim}_{\mathrm{AC},z}^{(n,1)}(t)\bigr)\,dW(t),
 \end{aligned}
\]
with initial condition
$\widetilde{\crlim}_{\mathrm{AC},z}^{(n,1)}(0)=z$.
Here $W$ is the $|\mathcal E|$-dimensional Brownian motion
used above. 
The omitted mutation covariance satisfies
\[
 \sup_{r\in\mathcal P(\mathbb S)}
 \left\|\bnG(r)-\crdiff_{\mathrm{AC}}(r)
                    \crdiff_{\mathrm{AC}}(r)^\top\right\|_{\mathrm{tr}}
 \leq\frac{2\Theta}{n}.
\]
We use the generator comparison argument in the proof of
Theorem~\ref{thm:critical_pairwise}, keeping this contribution
in the error. Let $\widetilde{\mathcal B}_{\mathrm{AC}}$
be the generator of the diffusion above. Following the proof of Theorem~\ref{thm:critical_pairwise}, we obtain that for
$f\in C_b^4(\R^K;\R)$, $x\in\mathbb S^n$, and $r=\nden(x)$,
\[
 \left|\bngen(f\circ\nden)(x)
       -\widetilde{\mathcal B}_{\mathrm{AC}}f(r)\right|
 \leq\frac{\Theta}{n}\|D^2f\|_{\mathrm{op},\infty}
      +\frac{C}{n^2}\|D^4f\|_{\mathrm{op},\infty}.
\]

Applying Lemma~\ref{lem:critical_backward_reg_PRM} and Dynkin's formula gives
\[
 \begin{aligned}
 \left|\EE[f(\bnmu_t)]
   -\EE\left[
     f\bigl(\widetilde{\crlim}_{\mathrm{AC},\mu_0^{(n)}}^{(n,1)}(t)\bigr)
   \right]\right|\leq C_T\left[
    \frac1n\sum_{j=1}^2\|D^jf\|_{\mathrm{op},\infty}
    +\frac1{n^2}\sum_{j=1}^4\|D^jf\|_{\mathrm{op},\infty}
    \right],
 \end{aligned}
\]
for all $0\leq t\leq T$.
Lemma~\ref{lem:regularization_dK} with $\ell=4$ gives,
for every $T>0$,
\[
 \sup_{0\leq t\leq T}
 \wass_1\bigl(
   \operatorname{Law}(\bnmu_t),
   \operatorname{Law}(
     \widetilde{\crlim}_{\mathrm{AC},\mu_0^{(n)}}^{(n,1)}(t))
 \bigr)
 \leq \widetilde C_{\mathrm{AC}}(T)n^{-1/2}.
\]
The constant depends only on $K,T,\Theta$ and
$\max_{i\to j\in\mathcal E}c_{i\to j}$, and is independent
of $n$ and the initial empirical state. Therefore, omitting the
mutation covariance does not change thr order of the quantitative bound.
\end{remark}

\subsubsection{Unbalanced conversion in the ODE regime}
\label{subsubsec:regulated_TK_ode_PRM}
Under the unbalanced conversion case, the conversion constants
$c_{i\to j}$ are not necessarily symmetric. The transition rate is
\begin{equation}
\label{eq:unbalanced_TK_individual_rates_PRM}
 H_{i,j}^{(n)}(r)=c_{i\to j}r_ir_j+\frac{\vartheta_jr_i}{n},
 \qquad i\ne j.
\end{equation}

We take $\alpha(n)=\beta(n)=1$, so we have $a_n=\sqrt n$.
We do not speed up time in this case, and we take
$\bnmu_t=\mu_t^{(n)}$.
The normalized drift and covariance are
\begin{align*}
 &\bar F^{(n)}(k,r)
 =r_k\sum_{j\ne k}r_j(c_{j\to k}r_j-c_{k\to j}r_k)
   +\frac{\vartheta_k\sum_{i=1}^K r_i^2-\Theta r_k^2}{n},\\
 &\bar G^{(n)}(r)
 =\sum_{i\ne j}r_i^2
   \left(c_{i\to j}r_j+\frac{\vartheta_j}{n}\right)
   \bar{\canvec}_{i,j}\bar{\canvec}_{i,j}^{\top}.
\end{align*}
Assumptions~\ref{assum:scaling-1} and~\ref{assum:ode-regime} hold.
Now we let $\rho^{(n)}$ solve \eqref{eq:ODE-approx}
with initial condition $\rho_0^{(n)}=\mu_0^{(n)}$.
On each fixed interval $[0,T]$,
Theorem~\ref{th:ODE-approx} bounds the $p$th moment of
$\sup_{t\leq T}\|\bnmu_t-\rho_t^{(n)}\|$.
Theorem~\ref{thm:kantorovich_PRM} gives a $C_Tn^{-1/6}$ bound in
$\wass_1$ between the law of fluctuation $\fluct_t^{(n)}$
and the time-$t$ law of the diffusion in \eqref{eq:Y_PRM_SDE}.
Taking $c_{\max}=\max_{i\to j\in\mathcal E}c_{i\to j}$,
Theorem~\ref{thm:concentration_PRM} applies with
$\|D\bar F\|_{\mathrm{op},\infty}\leq\sqrt{2K}\,c_{\max}+4\Theta$,
$\delta_T=\delta e^{-\|D\bar F\|_{\mathrm{op},\infty}T}$, and
\[
 \|\bar G\|_{\mathrm{tr},\infty}
 =2\sup_{n,r}\sum_{i\ne j}r_i^2
       \left(c_{i\to j}r_j+\frac{\vartheta_j}{n}\right)
 \leq 2(c_{\max}+\Theta).
\]
Thus \eqref{eq:concentration_main_PRM} gives the bound
$2e^{-C_{\mathrm{conc}}(T,\delta)n}$, where
\[
 C_{\mathrm{conc}}(T,\delta)=\frac{\delta_T^2}
 {2T\|\bar G\|_{\mathrm{tr},\infty}+(2\sqrt2/3)\delta_T}.
\]

For $K=2$, we identify species $1,2$ with types $0,1$.
In addition, we take
$c_{1\to2}=\kappa_1$, $c_{2\to1}=\kappa_2$,
$\vartheta_2=\iota_1$, and $\vartheta_1=\iota_2$.
\eqref{eq:unbalanced_TK_individual_rates_PRM} gives
\[
 \begin{aligned}
 H_{1,2}^{(n)}(r)
 =\kappa_1r_2(1-r_2)+\frac{\iota_1}{n}(1-r_2),\qquad
 H_{2,1}^{(n)}(r)
 =\kappa_2r_2(1-r_2)+\frac{\iota_2}{n}r_2.
 \end{aligned}
\]
The transition rates agree with those in
\cite[Section~2.1, equation~(4)]{KhudaBukhshXiang2026}.

\begin{appendix}
\section{Regularization results}
\label{app:regularization}
\subsection{Regularization for Wasserstein distance}
\label{Regularization_lemma_for_Wasserstein_distance}
We write $\mathcal P_1(\R^d)$ for the set of Borel probability
measures $\mu$ on $\R^d$ with finite first moment, that is,
$\int_{\R^d}\|x\|\,\mu(dx)<\infty$. The following lemma appeared in
\cite[Lemma~A.1]{AljovinJaraXiang2024} and is included here
for completeness.
\begin{lemma}[Regularization for Wasserstein distance]
\label{lem:regularization_dK}
Let $\mu,\nu\in\mathcal P_1(\R^d)$ and let $\ell\geq2$ be an
integer. Suppose that $b_1,\ldots,b_\ell\geq0$ satisfy
\begin{equation}
\label{eq:regularization_input}
 \left|\int f\,d\mu-\int f\,d\nu\right|
 \leq\sum_{j=1}^\ell b_j\|D^jf\|_{\mathrm{op},\infty},
 \qquad f\in C_b^\infty(\R^d;\R).
\end{equation}
Then there exists a finite constant $C_{\mathrm{reg}}(d,\ell)$,
depending only on $d$ and $\ell$, such that
\begin{equation}
\label{eq:regularization_conclusion}
 \wass_1(\mu,\nu)
 \leq C_{\mathrm{reg}}(d,\ell)
       \max_{1\leq j\leq\ell}b_j^{1/j}.
\end{equation}
\end{lemma}

\begin{proof}
We fix a $1$-Lipschitz function $f:\R^d\to\R$. We can assume $f(0)=0$, since subtracting $f(0)$ does not change its integral difference.
Then $|f(u)|\leq\|u\|$.

For $R>0$, we define
$$
 f^R(u)\dfeq\max\{-R,\min\{f(u),R\}\}.
$$
Then $f^R$ is bounded and $1$-Lipschitz. Moreover,
$|f^R(u)|\leq|f(u)|$ and $f^R(u)\to f(u)$ as $R\to\infty$.
We furthermore let
$$
 \gamma(u)\dfeq(2\pi)^{-d/2}e^{-\|u\|^2/2},
 \qquad
 \gamma_\theta(u)\dfeq\theta^{-d}\gamma(u/\theta),
 \qquad
 f^R_\theta\dfeq f^R*\gamma_\theta,
 \quad \theta>0.
$$
Since $f^R$ is bounded and every derivative of $\gamma_\theta$
is integrable, $f^R_\theta\in C_b^\infty(\R^d;\R)$.
The Lipschitz property of $f^R$ gives
\begin{align}
 \|f^R_\theta-f^R\|_\infty
 &\leq\int_{\R^d}\|v\|\gamma_\theta(v)\,dv
 =\theta\int_{\R^d}\|v\|\gamma(v)\,dv
 =C_d\theta,
 \label{eq:trunc_smoothing_error}
\end{align}
where
$C_d\dfeq\int_{\R^d}\|v\|\gamma(v)\,dv<\infty$.

For $j\geq1$ and $h_1,\ldots,h_j\in\R^d$ with
$\|h_i\|\leq1$, we notice that
 $\int_{\R^d}
 D^j\gamma_\theta(v)[h_1,\ldots,h_j]\,dv=0$.
It follows that
\begin{align*}
 D^jf^R_\theta(u)[h_1,\ldots,h_j]
 &=\int_{\R^d}
 \bigl(f^R(u-v)-f^R(u)\bigr)
 D^j\gamma_\theta(v)[h_1,\ldots,h_j]\,dv.
\end{align*}
Taking the supremum over $u$ gives
\begin{align}
 \|D^jf^R_\theta\|_{\mathrm{op},\infty}
 &\leq\int_{\R^d}
       \|v\|\|D^j\gamma_\theta(v)\|_{\mathrm{op}}\,dv
      =\theta^{1-j}\int_{\R^d}
       \|v\|\|D^j\gamma(v)\|_{\mathrm{op}}\,dv
 =C_{d,j}\theta^{1-j},
 \label{eq:derivative_mollified_bound}
\end{align}
where
$C_{d,j}\dfeq \int_{\R^d}\|v\|\|D^j\gamma(v)\|_{\mathrm{op}}\,dv<\infty$.
Thus the bounds in
\eqref{eq:trunc_smoothing_error} and
\eqref{eq:derivative_mollified_bound} are independent of
$R$ and $f$.

Applying \eqref{eq:regularization_input} to $f^R_\theta$
and using \eqref{eq:trunc_smoothing_error} and
\eqref{eq:derivative_mollified_bound}, we get
\begin{align}
 \left|\int f^R\,d\mu-\int f^R\,d\nu\right|
 &\leq 2\|f^R_\theta-f^R\|_\infty
       +\sum_{j=1}^\ell
        b_j\|D^jf^R_\theta\|_{\mathrm{op},\infty}
      \leq 2C_d\theta+
       \sum_{j=1}^\ell C_{d,j}b_j\theta^{1-j}.
 \label{eq:regularization_pre_optimization}
\end{align}
We take $A\dfeq\max_{1\leq j\leq\ell}b_j^{1/j}$.
If $A>0$, take $\theta=A$. It follows that
\begin{equation}
\label{eq:regularization_bound_truncated}
 \left|\int f^R\,d\mu-\int f^R\,d\nu\right|
 \leq\left(2C_d+\sum_{j=1}^\ell C_{d,j}\right)A.
\end{equation}
Since $|f^R(u)|\leq\|u\|$, we can use dominated convergence under
both measures by taking $R\to\infty$. Taking the supremum
over all $1$-Lipschitz functions with $f(0)=0$ proves
\eqref{eq:regularization_conclusion} with
$$
 C_{\mathrm{reg}}(d,\ell)
 \dfeq 2C_d+\sum_{j=1}^\ell C_{d,j}.
$$

If $A=0$, every $b_j$ is zero. Letting $\theta\downarrow0$
in \eqref{eq:regularization_pre_optimization} and then
$R\to\infty$ gives
$\int f\,d\mu=\int f\,d\nu$ for every $1$-Lipschitz
function $f$. Hence $\wass_1(\mu,\nu)=0$, completing
the proof.
\end{proof}

\subsection{Proof of the regularity lemma for the ODE regime}
\label{Regularity_for_ode_proof}
In this section, we prove Lemma \ref{lem:backward_reg_PRM}.

\begin{proof}[Proof of Lemma \ref{lem:backward_reg_PRM}]
We have $\Phi^{(n)}(t,s)$ is the fundamental matrix in \eqref{eq:Phi_PRM_def}. 
Since
$$
 \|\Phi^{(n)}(t,s)\|_{\mathrm{op}}
 \leq
 1+\int_s^t
 \|D\bnF(\rho_r^{(n)})\|_{\mathrm{op}}
 \|\Phi^{(n)}(r,s)\|_{\mathrm{op}}\,dr,
$$
Gronwall's inequality gives
\begin{align}\label{eq:fund-gronwall}
 \|\Phi^{(n)}(t,s)\|_{\mathrm{op}}
 &\leq
 \exp\left(
 \int_s^t\|D\bnF(\rho_r^{(n)})\|_{\mathrm{op}}\,dr
 \right)\notag\\
 &\leq
 \exp\left(\|D\bar F\|_{\mathrm{op},\infty}(t-s)\right).
\end{align}
By variation of constants, the solution to
\eqref{eq:Y_PRM_SDE} is given by
\begin{equation}\label{eq:linear_transition_ode}
 Y_{s,x}^{(n)}(t)
 =\Phi^{(n)}(t,s)x
 +\int_s^t\Phi^{(n)}(t,r)
 (\bnG)^{1/2}(\rho_r^{(n)})\,dW_r.
\end{equation}
Consequently,
$$
 g_{s,t}(x)
 =
 \EE\left[
 f\left(
 \Phi^{(n)}(t,s)x
 +\int_s^t\Phi^{(n)}(t,r)
 (\bnG)^{1/2}(\rho_r^{(n)})\,dW_r
 \right)
 \right].
$$
By differentiating under expectation, which is justified by boundedness of the derivatives of $f$, we obtain  (see \cite[Lemma~2.4, formula~(2.15)]{CerraiLunardi2021}), for
\(m=1,2,3\),
\begin{align}\label{eq:CL_derivative_formula_PRM}
 D^m g_{s,t}(x)[h_1,\ldots,h_m]= 
 \EE\left[
 D^m f(Y_{s,x}^{(n)}(t))
 [\Phi^{(n)}(t,s)h_1,\ldots,\Phi^{(n)}(t,s)h_m]
 \right].
\end{align}
Hence,
$$
 |D^m g_{s,t}(x)[h_1,\ldots,h_m]|
 \leq
 \|D^m f\|_{\mathrm{op},\infty}
 \|\Phi^{(n)}(t,s)\|_{\mathrm{op}}^m
 \prod_{i=1}^m\|h_i\|,
$$
which proves the first inequality in
\eqref{eq:backward_reg_PRM}. The second follows from
\eqref{eq:fund-gronwall}.
\end{proof}

\subsection{Proof of the regularity lemmas for the SDE regime}
\label{Regularity_for_Sde_proof}
In this section, we prove Lemma \ref{lem:critical_backward_reg_PRM} and Lemma \ref{lem:critical_WF_backward_regularity_PRM}.

\begin{proof}[Proof of Lemma \ref{lem:critical_backward_reg_PRM}]
Let $\ell=3$, or $\ell=4$ under the additional hypothesis,
and let $\psi\in C_b^\ell(\R^K;\R)$. Fix $c>0$.
The coefficients $\bnF$ and $c^{-1}\crdiff^{(n)}$ are globally
Lipschitz, have at most linear growth, and have bounded continuous
derivatives of orders $1,\ldots,\ell$. By
\cite[Chapter 5, Theorems 5.3--5.4 and the higher-derivative
remark on p.~123]{Friedman1975}, the spatial derivatives of
$\crlim_x^{(n,c)}$ exist in $L^2$ through order $\ell$ and
satisfy the corresponding variational equations.

Fix $n,x$ and directions $h_1,\ldots,h_\ell$ with $\|h_i\|\leq1$.
For each nonempty $I\subseteq\{1,\ldots,\ell\}$, put
$$
J_I(r)\dfeq
D_x^{|I|}\crlim_x^{(n,c)}(r)[h_i:i\in I],
$$
and let $\mathrm{Par}(I)$ denote the set of partitions of $I$.
The higher-order chain rule gives
\begin{align}\label{eq:critical_variational_PRM}
\begin{aligned}
dJ_I(r)
={}&\sum_{\mathscr P\in\mathrm{Par}(I)}
D^{|\mathscr P|}\bnF(\crlim_x^{(n,c)}(r))
[J_B(r):B\in\mathscr P]\,dr\\
&+c^{-1}\sum_{\mathscr P\in\mathrm{Par}(I)}
D^{|\mathscr P|}\crdiff^{(n)}(\crlim_x^{(n,c)}(r))
[J_B(r):B\in\mathscr P]\,dW_r,
\end{aligned}
\end{align}
where $J_{\{i\}}(0)=h_i$ and $J_I(0)=0$ for $|I|\geq2$.
We use the continuous versions supplied by these equations.

We prove moment bounds by induction on $m=|I|$, simultaneously
for every $p\geq2$. The partition $\{I\}$ gives the terms
linear in $J_I$; all other partitions contain only proper
subsets of $I$. Set
$$
\tau_N\dfeq
\inf\left\{r\geq0:
\max_{\varnothing\neq B\subseteq I}\|J_B(r)\|\geq N
\right\}\wedge T,
\qquad N\geq2.
$$
It\^o's formula for $\|J_I(r\wedge\tau_N)\|^p$, Young's
inequality, and
$\trace(\sigma\sigma^{\mathsf T})\leq K\|\sigma\|_{\mathrm{op}}^2$
give
\begin{align}\label{eq:critical_flow_moment_recursion_PRM}
\begin{aligned}
\EE\|J_I(r\wedge\tau_N)\|^p
\leq{}&1+\Gamma_{m,p}(c)\int_0^r
\EE\|J_I(u\wedge\tau_N)\|^p\,du\\
&+\Gamma_{m,p}(c)
\sum_{\substack{\mathscr P\in\mathrm{Par}(I)\\|\mathscr P|\geq2}}
\int_0^r
\EE\left[\prod_{B\in\mathscr P}\|J_B(u)\|^p\right]\,du,
\end{aligned}
\end{align}
where, for a sufficiently large constant $\gamma_{m,p,K}$
depending only on $m,p,K$, we may take
$$
\Gamma_{m,p}(c)
\dfeq\gamma_{m,p,K}\left[
\max_{1\leq j\leq m}\|D^j\bar F\|_{\mathrm{op},\infty}
+c^{-2}\max_{1\leq j\leq m}
\|D^j\Sigma\|_{\mathrm{op},\infty}^{\,2}
\right].
$$
Indeed, the drift terms are linear in the bounds on the
derivatives of $\bnF$, while the quadratic-variation terms
are quadratic in those of $c^{-1}\crdiff^{(n)}$.
The stopped stochastic integral has mean zero.

For $m=1$, the sum in
\eqref{eq:critical_flow_moment_recursion_PRM} is empty.
For $m\geq2$, H\"older's inequality and the induction hypothesis give
$$
\EE\left[\prod_{B\in\mathscr P}\|J_B(u)\|^p\right]
\leq
\prod_{B\in\mathscr P}
\left(\EE\|J_B(u)\|^{p|\mathscr P|}\right)^{1/|\mathscr P|}
\leq
\prod_{B\in\mathscr P}
C_{|B|,p|\mathscr P|,c}(T)^{1/|\mathscr P|}.
$$
By Gronwall's inequality, we have a bound independent of
$N,n,x$ and the unit-bounded directions. By passing $N\to\infty$
and using Fatou's lemma, we obtain
\begin{equation}\label{eq:critical_flow_derivative_bound_PRM}
\sup_{\substack{n\geq1,\ 0\leq r\leq T,\ x\in\R^K\\
\|h_1\|,\ldots,\|h_m\|\leq1}}
\EE\left\|
D_x^m\crlim_x^{(n,c)}(r)[h_1,\ldots,h_m]
\right\|^p
\leq C_{m,p,c}(T),
\qquad 1\leq m\leq\ell.
\end{equation}
Moreover, $\Gamma_{m,p}(c)$ is nonincreasing in $c$.
Choosing at each induction step the bound supplied by
Gronwall's inequality therefore makes $C_{m,p,c}(T)$
nonincreasing in $c$: the recursion uses only sums, products,
positive powers, and exponentials of nonnegative quantities.
For $1\leq p<2$, take
$C_{m,p,c}(T)=C_{m,2,c}(T)^{p/2}$.

By \cite[Chapter 5, Theorem 5.5 and the higher-derivative
remark on p.~123]{Friedman1975}, differentiation under the
expectation is valid. For $I=\{1,\ldots,m\}$,
\begin{equation}\label{eq:critical_backward_derivative_formula_PRM}
D^m\bar g_{s,t}(x)[h_1,\ldots,h_m]
=
\EE\left[
\sum_{\mathscr P\in\mathrm{Par}(I)}
D^{|\mathscr P|}\psi(\crlim_x^{(n,c)}(t-s))
[J_B(t-s):B\in\mathscr P]
\right].
\end{equation}
H\"older's inequality and
\eqref{eq:critical_flow_derivative_bound_PRM} yield
\begin{align}\label{eq:critical_backward_derivative_bound_PRM}
|D^m\bar g_{s,t}(x)[h_1,\ldots,h_m]|
&\leq
\sum_{\mathscr P\in\mathrm{Par}(I)}
\|D^{|\mathscr P|}\psi\|_{\mathrm{op},\infty}
\prod_{B\in\mathscr P}
C_{|B|,|\mathscr P|,c}(T)^{1/|\mathscr P|}\notag\\
&\leq
C_{\mathrm{back},m}(c,T)
\sum_{j=1}^m\|D^j\psi\|_{\mathrm{op},\infty},
\end{align}
where we may take
$$
C_{\mathrm{back},m}(c,T)
\dfeq
\sum_{\mathscr P\in\mathrm{Par}(\{1,\ldots,m\})}
\prod_{B\in\mathscr P}
C_{|B|,|\mathscr P|,c}(T)^{1/|\mathscr P|}.
$$
This constant is finite and nonincreasing in $c$.
Taking the supremum over $x$ and the unit-bounded directions
proves the derivative estimate. Also,
\begin{equation}\label{eq:critical_backward_zero_bound_PRM}
\|\bar g_{s,t}\|_\infty\leq\|\psi\|_\infty.
\end{equation}
Continuity of the derivative processes and their
higher-moment bounds imply, by uniform integrability,
continuity in $(s,t)$ of the derivatives in
\eqref{eq:critical_backward_derivative_formula_PRM}
for each fixed $x$, including at $s=t$.
At $s=t$, only the partition into singleton blocks contributes. This  gives $D^m\bar g_{t,t}(x)=D^m\psi(x)$.

Finally, $\bnF$, $c^{-1}\crdiff^{(n)}$, and $\psi$ satisfy
the hypotheses of \cite[Chapter 5, Theorem 6.1]{Friedman1975}.
Therefore, $\bar g_{s,t}$ satisfies the classical backward equation
\eqref{eq:critical_bke_PRM}, with $\bar g_{t,t}=\psi$.

Recall that $a_n>0$ and $a_n\to a_\infty\in(0,\infty)$.
We take $a_*\dfeq\inf_{n\geq1}a_n>0$. By the monotone choice of constants above, we have that
$$
\sup_{n\geq1}C_{\mathrm{back},m}(a_n,T)
\leq C_{\mathrm{back},m}(a_*,T)<\infty,
\qquad 1\leq m\leq\ell.
$$
\end{proof}

Now we prove Lemma \ref{lem:critical_WF_backward_regularity_PRM}.
\begin{proof}[Proof of Lemma \ref{lem:critical_WF_backward_regularity_PRM}]
We fix $n,c$, and write $b=\mathfrak p\circ\bnF\circ\mathfrak r$. We take
\[
 Ah(x)=\sum_{i,j=1}^d x_i(\delta_{ij}-x_j)\partial_{ij}h(x),
 \qquad L=\mathcal B_{\mathrm{WF},c}^{(n)}=c^{-2}A+b\cdot D.
\]
We define 
$p(x):=\det(x_i(\delta_{ij}-x_j))_{i,j=1}^d
=\mathfrak r_K(x)\prod_{j=1}^d x_j$.
Since $\partial_i\mathfrak r_K(x)=-1$, the product rule yields
\[
 b(x)\cdot\nabla p(x)=
 \begin{cases}
  b_i(x)\,\mathfrak r_K(x)\displaystyle\prod_{j\ne i}x_j,
       & x_i=0,\\[3pt]
  -\bigl(\displaystyle\sum_{i=1}^d b_i(x)\bigr)
       \displaystyle\prod_{j=1}^d x_j,
       & \mathfrak r_K(x)=0.
 \end{cases}
\]
Both cases have nonnegative expression.
Thus $c^2b/2$ satisfies condition~(6.2) of \cite{Campiti2002}.

Let $P^{(n,c)}=(P_u^{(n,c)})_{u\geq0}$ denote the transition
semigroup of the first $d$ coordinates of the Wright--Fisher diffusion \eqref{eq:WF_diffusion_SDE_PRM}. It follows that
\[
 (P_u^{(n,c)}h)(x)
 \dfeq
 \EE\left[
 h\left(
 \mathfrak p\bigl(
 \crlim_{\mathrm{WF},\mathfrak r(x)}^{(n,c)}(u)
 \bigr)
 \right)
 \right],
 \qquad h\in C(E;\R),\quad u\geq0,\quad x\in E.
\]
In particular, $\bar g_{s,t}=P_{t-s}^{(n,c)}\psi$.
Applying Theorem~6.1 of \cite{Campiti2002} to $c^2L/2$ and rescaling
time by $2c^{-2}$, we obtain that the closure of $L$ on $C^2(E)$
generates $P^{(n,c)}$ and that $P^{(n,c)}$ preserves $C^m(E)$
for $1\leq m\leq4$. In particular,
$\partial_uP_u^{(n,c)}\psi=LP_u^{(n,c)}\psi$ on $E$. Taking $u=t-s$
gives the backward equation with terminal value $\psi$.
For the comparison diffusion
$X_u=\mathfrak p(\crlim_{\mathrm{WF},\mathfrak r(x)}^{(n,c)}(u))$
started at $x\in E$, we have
$\bar g_{s,t}(x)=\mathbb E_x\!\left[\psi(X_{t-s})\right]$. It follows that
\[
 |\bar g_{s,t}(x)|
 \leq \mathbb E_x\!\left[|\psi(X_{t-s})|\right]
 =\bigl(P_{t-s}^{(n,c)}|\psi|\bigr)(x)
 \leq \|\psi\|_\infty.
\]
We have
$\|\bar g_{s,t}\|_\infty\leq\|\psi\|_\infty$ by taking the supremum over $x\in E$.

We let $Q_u$
be the drift-free semigroup with generator $A$. We also write
$H_m(h)=\max_{1\leq|\nu|\leq m}\|\partial^\nu h\|_\infty$.
To apply \cite[Lemma~3.1]{ChenXiongXuZheng2026}, we take $K$ alleles and take all mutation parameters to zero. In the notation of
equations~(2.3)--(2.7), the generator then
reduces to $A/(4N)$. The decay exponent for a mixed derivative of order $k$ is $k(k-1)/(4N)$. Applying estimate~(3.4) to $h$
and $-h$, at time $4Nu$ gives, for $h\in C^4(E)$,
\[
 \|\partial^\nu Q_u h\|_\infty
 \leq e^{-k(k-1)u}\|\partial^\nu h\|_\infty,
 \qquad k=|\nu|,\quad 1\leq k\leq4.
\]
Hence $H_m(Q_{c^{-2}u}h)\leq H_m(h)$ for every $u\geq0$.

For each initial state $x\in E$, we
define $\phi_u^{(n)}(x)$ as the solution at time $u\geq0$ of
\[
 \partial_u\phi_u^{(n)}(x)=b(\phi_u^{(n)}(x)),
 \qquad \phi_0^{(n)}(x)=x.
\]
For fixed $u$, its Jacobian with
respect to the initial state is the $d\times d$ matrix
\[
 J_u(x)=D_x\phi_u^{(n)}(x)
 =\bigl(\partial_{x_j}\phi_{u,i}^{(n)}(x)\bigr)_{1\leq i,j\leq d}.
\]
By the $C^6$ regularity of $b$, these derivatives exist  and can be extended continuously to the boundary of $E$. Differentiating with respect to $x$ gives
\[
 \partial_u J_u(x)=Db(\phi_u^{(n)}(x))J_u(x),
 \qquad J_0(x)=I_d.
\]
By Gronwall's inequality, we obtain
$\|J_u\|_{\mathrm{op},\infty}
\leq e^{u\|Db\|_{\mathrm{op},\infty}}$.
By integrating the equation for $J_u$, it follows that
$\|J_u-I_d\|_{\mathrm{op},\infty}
\leq e^{u\|Db\|_{\mathrm{op},\infty}}-1\leq C_1u$
for $0\leq u\leq1$. At second order, we have that
\[
 \frac{d}{du}D^2\phi_u^{(n)}[v,w]
 =Db(\phi_u^{(n)})D^2\phi_u^{(n)}[v,w]
  +D^2b(\phi_u^{(n)})[J_uv,J_uw].
\]
For unit $v,w$, the last term is uniformly bounded on $[0,1]$. For higher orders $2\leq j\leq m$, we differentiate $j$ times
with respect to the initial point $x$. On the right-hand side,
$D_x^j\phi_u^{(n)}(x)$ appears only once, multiplied by the
Jacobian matrix $Db(\phi_u^{(n)}(x))$. In every other term, the flow is differentiated at most $j-1$ times. So its derivatives
are already bounded by induction. The derivatives of $b$ in
those terms have orders between $2$ and $j$ and are uniformly
bounded by hypothesis. Using the fact that $D^j\phi_0^{(n)}=0$, and by induction and
Gronwall's inequality, we have 
\[
 \|D\phi_u^{(n)}-I_d\|_{\mathrm{op},\infty}
 +\sum_{j=2}^m\|D^j\phi_u^{(n)}\|_{\mathrm{op},\infty}
 \leq C_m u,\qquad 0\leq u\leq1.
\]
$C_m$ depends only on $d,m$ and the uniform bounds
on the drift derivatives of orders $1,\ldots,m$.

To estimate $h\circ\phi_u^{(n)}$, we write
$D\phi_u^{(n)}=I_d+(D\phi_u^{(n)}-I_d)$.
At the first order, by using the chain rule, we obtain
\[
 \partial_i(h\circ\phi_u^{(n)})
 =(\partial_i h)\circ\phi_u^{(n)}
  +\sum_{j=1}^d \bigl((\partial_jh)\circ\phi_u^{(n)}\bigr)
                (\partial_i\phi_{u,j}^{(n)}-\delta_{ij}).
\]
At any order $1\leq|\nu|\leq m$, by taking the identity matrix in every first-derivative factor, we have that
$(\partial^\nu h)\circ\phi_u^{(n)}$
with coefficient one. Each remaining term contains either
$D\phi_u^{(n)}-I_d$ or a derivative $D^j\phi_u^{(n)}$ with $j\geq2$. Therefore, they are bounded by $C_m u H_m(h)$.
Since $\phi_u^{(n)}(E)\subseteq E$,
taking the maximum over $\nu$ yields
\[
 H_m(h\circ\phi_u^{(n)})\leq(1+C_m u)H_m(h),
 \qquad 0\leq u\leq1.
\]

The drift-free diffusion with time rescaled by $c^{-2}$ defines
a strongly continuous contraction semigroup on $C(E)$.
Composition with the deterministic drift flow defines another such semigroup.
Their generators act as $c^{-2}A$ and $b\cdot D$ on $C^2(E)$,
respectively. It follows that their sum agrees with $L$.
By \cite[Theorem~6.1]{Campiti2002}, the closure of $L$ on
$C^2(E)$ generates $P^{(n,c)}$. For $u>0$, the product
formula of \cite{Trotter1959} gives the uniform limit on $E$
\[
 P_u^{(n,c)}\psi
 =\lim_{N\to\infty} Q_{c^{-2}u/N}\!\Bigl(
   \bigl[\cdots Q_{c^{-2}u/N}
     (\psi\circ\phi_{u/N}^{(n)})\cdots\bigr]
   \circ\phi_{u/N}^{(n)}\Bigr).
\]
Each of the $N$ steps composes the current test function with
$\phi_{u/N}^{(n)}$ and then applies $Q_{c^{-2}u/N}$.
For $u/N\leq1$, the preceding estimates give
\[
 H_m\!\left(Q_{c^{-2}u/N}
     \bigl(h\circ\phi_{u/N}^{(n)}\bigr)\right)
 \leq H_m\bigl(h\circ\phi_{u/N}^{(n)}\bigr)
 \leq (1+C_m u/N)H_m(h).
\]
After $N$ steps the seminorm is therefore bounded by
$(1+C_m u/N)^N H_m(\psi)\leq e^{C_m u}H_m(\psi)$.

Since $P_u^{(n,c)}\psi\in C^4(E)$, the bound extends to the
boundary, giving
$H_m(P_u^{(n,c)}\psi)\leq e^{C_m u}H_m(\psi)$.
This also holds at $u=0$ since $P_0^{(n,c)}\psi=\psi$. Taking
$u=t-s\leq T$ and operator norms proves \eqref{eq:WF_backward_regularity_PRM}.

Finally, for the scalar drift $A_n(1-x)-B_nx$, the operator $c^2L$
is the Wright--Fisher generator in
\cite[Lemma~3.3]{AljovinJaraXiang2024} with mutation parameters
$c^2A_n$ and $c^2B_n$. At time $c^{-2}(t-s)$, we have that
\[
 \|\partial_x^m\bar g_{s,t}\|_\infty
 \leq e^{-m[A_n+B_n+c^{-2}(m-1)](t-s)}
       \|\psi^{(m)}\|_\infty,
\]
which proves \eqref{eq:WF_neutral_derivative_decay_PRM}.
\end{proof}

\section{Auxiliary Lemma}
\label{app:Auxiliary Lemma}
\begin{lemma}[Burkholder--Rosenthal type Inequality]
\label{lem:burkholder_rosenthal_PRM}
Let $M$ be an $\R^K$-valued c\`adl\`ag locally
square-integrable martingale with $M(0)=0$.
Then for any $p \geq 2$, there exists a constant $C_{\mathrm{BR}}(p)>0$ depending only on $p$, such that
\[
\mathbb E\!\left(\sup_{s\le t}\|M(s)\|^p\right)
\le C_{\mathrm{BR}}(p)\lf(\mathbb E\!\left(
 \bigl(\operatorname{tr}\langle M\rangle_t\bigr)^{p/2}\right)
+\mathbb E\!\left(
 \sup_{0<s\le t}\|\Delta M(s)\|^p\right)\ri),
\]
\end{lemma}

A proof can be found in \cite[Theorem~3.6]{MarinelliRockner2014}. 
For the scalar version, see
\cite[Lemma~2.1, equation~(2.7)]{DzhaparidzeValkeila1990}.

\end{appendix}

\bibliographystyle{spr-ims-number}
\bibliography{references}
\end{document}